\documentclass[a4paper,twoside]{article}

\makeatletter
\renewcommand{\@seccntformat}[1]{\csname the#1\endcsname}
\renewcommand\normalsize{%
   \@setfontsize\normalsize\@xpt\@xiipt
   \abovedisplayskip 6\p@ \@plus1\p@
   \abovedisplayshortskip 6\p@ \@plus1\p@
   \belowdisplayshortskip 6\p@ \@plus1\p@
   \belowdisplayskip \abovedisplayskip
   \let\@listi\@listI}
\normalsize
\renewcommand\ps@headings{%
  \let\@oddfoot\@empty\let\@evenfoot\@empty
  \let\@oddhead\@empty\let\@evenhead\@empty
  \let\@mkboth\markboth
  \def\sectionmark##1{%
    \markboth {\uppercase{\ifnum \c@secnumdepth >\z@
        \thesection.\relax\fi
        ##1}}{}}%
  \def\subsectionmark##1{%
    \markright {\ifnum \c@secnumdepth >\@ne
        \thesubsection\relax \fi
        ##1}}}
\renewcommand\ps@myheadings{%
    \let\@oddfoot\@empty\let\@evenfoot\@empty
    \def\@evenhead{\thepage\hfil\small\leftmark\hfil}%
    \def\@oddhead{\hfil{\small\rightmark}\hfil\thepage}%
    \let\@mkboth\@gobbletwo
    \let\sectionmark\@gobble
    \let\subsectionmark\@gobble
    }
\def\myfnsymbol#1{\expandafter\@myfnsymbol\csname c@#1\endcsname}
\def\@myfnsymbol#1{\ensuremath{\ifcase#1\or *\or \ddagger\or **\or
   \mathsection\or \mathparagraph\or \|\or **\or \dagger\dagger
   \or \ddagger\ddagger \else\@ctrerr\fi}}
\renewcommand\maketitle{\par
  \begingroup
    \renewcommand\thefootnote{\myfnsymbol{footnote}}%
    \def\@makefnmark{\hbox to\z@{$\m@th^{\@thefnmark}$\hss}}%
    \long\def\@makefntext##1{\parindent 1em\noindent
            \hbox to1.8em{\hss$\m@th^{\@thefnmark}$}##1}%
    \if@twocolumn
      \ifnum \col@number=\@ne
        \@maketitle
      \else
        \twocolumn[\@maketitle]%
      \fi
    \else
      \newpage
      \global\@topnum\z@   
      \@maketitle
    \fi
    \thispagestyle{headings}
    \@thanks
  \endgroup
  \setcounter{footnote}{0}%
  \let\thanks\relax
  \let\maketitle\relax\let\@maketitle\relax
  \gdef\@thanks{}\gdef\@author{}\gdef\@title{}}
\renewcommand\@maketitle{%
  \newpage
  \begin{center}%
    {\LARGE\bfseries \@title \par}%
    \vskip 24\p@%
        {\Large\itshape
      \lineskip .5em%
      \begin{tabular}[t]{c}%
        \@author
      \end{tabular}\par}%
  \end{center}%
  \par
  \vskip 60\p@}
\def\@sect#1#2#3#4#5#6[#7]#8{\ifnum #2>\c@secnumdepth
     \let\@svsec\@empty\else
     \refstepcounter{#1}%
     \let\@@protect\protect
     \def\protect{\noexpand\protect\noexpand}%
     \edef\@svsec{\@seccntformat{#1}}%
     \let\protect\@@protect\fi
     \@tempskipa #5\relax
      \ifdim \@tempskipa>\z@
      {#6\relax{\interlinepenalty \@M \@svsec #8\par}}%
               \csname #1mark\endcsname{#7}\addcontentsline
         {toc}{#1}{\ifnum #2>\c@secnumdepth \else
                      \protect\numberline{\csname the#1\endcsname}\fi
                    #7}\else
        \def\@svsechd{#6%
\hskip #3\relax  
                   \@svsec #8\csname #1mark\endcsname
                      {#7}\addcontentsline
                           {toc}{#1}{\ifnum #2>\c@secnumdepth \else
                           \protect\numberline{\csname the#1\endcsname}%
                                     \fi
                       #7}}\fi
     \@xsect{#5}}
\def\@ssect#1#2#3#4#5{\@tempskipa #3\relax
   \ifdim \@tempskipa>\z@
     {#4%
        {\interlinepenalty \@M #5\par}%
     }
   \else \def\@svsechd{#4%
\hskip #1
    \relax
                        #5}\fi
    \@xsect{#3}}

\renewcommand\section{\@startsection {section}{1}{\z@}%
                                  {-36\p@ \@plus -1\p@ \@minus -4\p@}%
                                  {12\p@ \@plus 1\p@}%
{\reset@font\Large\bfseries\centering}}
\renewcommand\subsection{\@startsection{subsection}{2}{\z@}%
                                     {-24\p@ \@plus -1\p@ \@minus -4\p@}%
                                     {8\p@ \@plus 1\p@}
                                     {\reset@font\large\bfseries\centering}}

\renewcommand\subsubsection{\@startsection{subsubsection}{3}{\z@}%
                                     {3.25ex \@plus 1ex \@minus .2ex}%
                                     {-1.5ex \@plus .2ex}%
                                     {\reset@font\normalsize\bfseries}}

\renewenvironment{abstract}{\noindent\small{\bfseries\abstractname.}%
        }%
        {\vskip 11\p@}
\newenvironment{classification}{\noindent\small 2000 Mathematics 
Subject Classification:}{\vskip 12\p@}
\def\@lbibitem[#1]#2{\item[\hfill\@biblabel{#1}]\if@filesw
      {\let\protect\noexpand
       \immediate
       \write\@auxout{\string\bibcite{#2}{#1}}}\fi\ignorespaces}

\renewenvironment{thebibliography}[1]
     {\section*{\reset@font\fontsize{11.6}{13.6pt}\bfseries\refname
        \@mkboth{\uppercase{\refname}}{\uppercase{\refname}}}%
      \list{\@biblabel{\arabic{enumiv}}}%
           {\settowidth\labelwidth{\@biblabel{#1}}%
            \setlength{\leftmargin}{\labelwidth}
            \labelsep2mm
            \itemsep4pt
        \topsep\z@
            \parsep\z@
            \advance\leftmargin\labelsep
                        \usecounter{enumiv}%
            \let\p@enumiv\@empty
            }%
            \sloppy\clubpenalty4000\widowpenalty4000%
      \sfcode`\.=\@m
      \small}
     {\def\@noitemerr
       {\@latex@warning{Empty `thebibliography' environment}}%
      \endlist}

\def\@seccntformat#1{\csname the#1\endcsname.\hskip 0.25em}
\newskip\aline \newskip\halfaline
\def\skipaline{\vskip\aline}

\def\qedbox{$\rlap{$\sqcap$}\sqcup$}

\def\Proof{\ifdim\lastskip<\aline\removelastskip\skipaline\fi
\noindent\it Proof. \rm}

\RequirePackage{amsmath}
\RequirePackage{amssymb}
\RequirePackage{ifthen}
\RequirePackage{xypic}

\def\FMithmInfo{1995/11/23 v2.2c Theorem extension package (FMi)}
\@ifundefined{theorem@style}{}{\endinput}

\gdef\theoremstyle#1{%
   \@ifundefined{th@#1}{\@warning
          {Unknown theoremstyle `#1'. Using `plain'}%
          \theorem@style{plain}}%
      {\theorem@style{#1}}%
      \begingroup
        \csname th@\the\theorem@style \endcsname
      \endgroup}
\global\let\@begintheorem\relax
\global\let\@opargbegintheorem\relax
\newtoks\theorem@style
\global\theorem@style{plain}
\gdef\theorembodyfont#1{%
   \def\@tempa{#1}%
   \ifx\@tempa\@empty
    \theorem@bodyfont{}%
   \else
    \theorem@bodyfont{\reset@font#1}%
   \fi
   }
\newtoks\theorem@bodyfont
\global\theorem@bodyfont{}
\gdef\theoremheaderfont#1{\gdef\theorem@headerfont{#1}%
       \gdef\theoremheaderfont##1{%
        \typeout{\string\theoremheaderfont\space should be used
                 only once.}}}
\ifx\upshape\undefined
\gdef\theorem@headerfont{\bfseries}
\else \gdef\theorem@headerfont{\normalfont\bfseries}\fi
\gdef\th@plain{\@input@{thp.sty}}
\gdef\th@break{\@input@{thb.sty}}
\gdef\th@marginbreak{\@input@{thmb.sty}}
\gdef\th@changebreak{\@input@{thcb.sty}}
\gdef\th@change{\@input@{thc.sty}}
\gdef\th@changep{\@input@{thcp.sty}}
\gdef\th@changepd{\@input@{thcpd.sty}}
\gdef\th@margin{\@input@{thm.sty}}
\gdef\@xnthm#1#2[#3]{\expandafter\@ifdefinable\csname #1\endcsname
   {%
    \@definecounter{#1}\@newctr{#1}[#3]%
    \expandafter\xdef\csname the#1\endcsname
      {\expandafter \noexpand \csname the#3\endcsname
       \@thmcountersep \@thmcounter{#1}}%
    \def\@tempa{\global\@namedef{#1}}%
    \expandafter \@tempa \expandafter{%
      \csname th@\the \theorem@style
            \expandafter \endcsname \the \theorem@bodyfont
     \@thm{#1}{#2}}%
    \global \expandafter \let \csname end#1\endcsname \@endtheorem
   }}
\gdef\@ynthm#1#2{\expandafter\@ifdefinable\csname #1\endcsname
   {\@definecounter{#1}%
    \expandafter\xdef\csname the#1\endcsname{\@thmcounter{#1}}%
    \def\@tempa{\global\@namedef{#1}}\expandafter \@tempa
     \expandafter{\csname th@\the \theorem@style \expandafter
     \endcsname \the\theorem@bodyfont \@thm{#1}{#2}}%
    \global \expandafter \let \csname end#1\endcsname \@endtheorem}}
\gdef\@othm#1[#2]#3{%
  \expandafter\ifx\csname c@#2\endcsname\relax
   \@nocounterr{#2}%
  \else
   \expandafter\@ifdefinable\csname #1\endcsname
   {\expandafter \xdef \csname the#1\endcsname
     {\expandafter \noexpand \csname the#2\endcsname}%
    \def\@tempa{\global\@namedef{#1}}\expandafter \@tempa
     \expandafter{\csname th@\the \theorem@style \expandafter
     \endcsname \the\theorem@bodyfont \@thm{#2}{#3}}%
    \global \expandafter \let \csname end#1\endcsname \@endtheorem}%
  \fi}
\gdef\@thm#1#2{\refstepcounter{#1}%
   \trivlist
   \@topsep \theorempreskipamount               
   \@topsepadd \theorempostskipamount           
   \@ifnextchar [%
   {\@ythm{#1}{#2}}%
   {\@begintheorem{#2}{\csname the#1\endcsname}\ignorespaces}}
\global\let\@xthm\relax
\newskip\theorempreskipamount
\newskip\theorempostskipamount
\global
\global
\global\let\@endtheorem=\endtrivlist
\@onlypreamble\@xnthm
\@onlypreamble\@ynthm
\@onlypreamble\@othm
\@onlypreamble\newtheorem
\@onlypreamble\theoremstyle
\@onlypreamble\theorembodyfont
\@onlypreamble\theoremheaderfont
\theoremstyle{plain}

\begingroup \makeatletter
\@ifundefined{theorem@style}{\input{theorem.sty}}{}
             [\FMithmInfo]
\gdef\th@plain{\normalfont\itshape
  \def\@begintheorem##1##2{%
        \item[\hskip\labelsep \theorem@headerfont ##1\ ##2.]}%
\def\@opargbegintheorem##1##2##3{%
   \item[\hskip\labelsep \theorem@headerfont ##1\ ##2\ \rm(##3).]}}
\endgroup

\def\@addpunct#1{\ifnum\spacefactor>\@m \else#1\fi}

\newenvironment{proof}[1][\proofname]{\par
  \normalfont
  \topsep6\p@\@plus6\p@ \trivlist
  \item[\hskip\labelsep\itshape
    #1\@addpunct{.}]\ignorespaces
}{\proof@ending\endtrivlist\par}

\newcommand{\proofname}{Proof}
\newcommand{\proof@ending}{\hfill\qedbox} \newcommand{\BoxedProofs}{\renewcommand{\proof@ending}{\hfill\ensuremath{\Box}}}
\newcommand{\NonBoxedProofs}{\renewcommand{\proof@ending}{\hfill\qedbox}}

\newskip\aline \newskip\halfaline
\def\skipaline{\vskip\aline}

\newcommand{\@lstlabel}{}
\newcommand{\lstlabel}[1]{\renewcommand{\@lstlabel}{#1}}
\newcommand{\@lsttemplate}{}
\newcommand{\lsttemplate}[1]{\renewcommand{\@lsttemplate}{#1}}
{\renewcommand{\theenumi}{\@lstlabel}%
  \begin{list}%
    {\@lstlabel}%
    {\usecounter{enumi}%
      \settowidth{\labelwidth}{\@lsttemplate}%
      \setlength{\leftmargin}{\labelwidth}%
      \setlength{\labelsep}{0pt}%
      \setlength{\topsep}{\medskipamount}%
      \setlength{\parsep}{0pt}%
      \setlength{\itemsep}{0pt}%
      \setlength{\itemindent}{0pt}%
      \setlength{\listparindent}{\parindent}}}
  {\end{list}}
\newcommand{\deflststyle}[2]{\expandafter\newcommand\expandafter{\csname @@#1@@\endcsname}{#2}}
\newcommand{\lststyle}[1]{\csname @@#1@@\endcsname}
\deflststyle{}{
  \lsttemplate{}
  \lstlabel{}}
\deflststyle{ }{
  \lsttemplate{\hspace{\parindent}}
  \lstlabel{}}
\deflststyle{--}{
  \lsttemplate{\hspace{\parindent}}
  \lstlabel{\,--\hfill}}
\deflststyle{a.}{
  \lsttemplate{\textnormal{\,b.\ \ }}
  \lstlabel{\textnormal{\,\alph{enumi}.\hfill}}}
\deflststyle{a)}{
  \lsttemplate{\textnormal{\,b). }}
  \lstlabel{\textnormal{\,\alph{enumi})\hfill}}}
\deflststyle{(a)}{
  \lsttemplate{\textnormal{(b)\ \ }}
  \lstlabel{\textnormal{(\alph{enumi})\hfill\ \ }}}
\deflststyle{1.}{
  \lsttemplate{\textnormal{\,0.\ \ }}
  \lstlabel{\textnormal{\hfill\arabic{enumi}.\ \ }}}
\deflststyle{1)}{
  \lsttemplate{\textnormal{\,0)\ \ }}
  \lstlabel{\textnormal{\hfill\arabic{enumi})\ \ }}}
\deflststyle{(1)}{
  \lsttemplate{\textnormal{(0)\ \ }}
  \lstlabel{\textnormal{\hfill(\arabic{enumi})\ \ }}}
\deflststyle{i.}{
  \lsttemplate{\textnormal{\,iii.\ \ }}
  \lstlabel{\textnormal{\hfill\roman{enumi}.\ \ }}}
\deflststyle{i)}{
  \lsttemplate{\textnormal{\,iii)\ \ }}
  \lstlabel{\textnormal{\hfill\roman{enumi})\ \ }}}
\deflststyle{(i)}{
  \lsttemplate{\textnormal{(iii)\ \ }}
  \lstlabel{\textnormal{\hfill(\roman{enumi})\ \ }}}
\lststyle{(a)} 

\newenvironment{block}{\begin{trivlist}\item{}}{\end{trivlist}}

\reversemarginpar

\let\s@bsection=\subsection
\newcommand{\c@mda}[2][]{\s@bsection[#1]{#2.}}
\newcommand{\c@mdb}[1]{\s@bsection*{#1.}}
\renewcommand{\subsection}{\secdef\c@mda\c@mdb}

\let\oldcite=\cite
\renewcommand{\cite}[2][no@ption]
{\ifthenelse{\equal{#1}{no@ption}}{\textnormal{\oldcite{#2}}}{\oldcite[#1]{#2}}}
\let\oldref=\ref
\renewcommand{\ref}[1]{\textnormal{\oldref{#1}}}
\let\oldpageref=\ref
\renewcommand{\pageref}[1]{\textnormal{\oldpageref{#1}}}

\let\@contact=\empty
\newcommand{\contact}[2][]{
  \expandafter\gdef\expandafter\@contact\expandafter{%
    \@contact
  \small
  \begin{block}#2\\[.5em]
  \ifthenelse{\equal{#1}{}}{}{Email: #1}\end{block}
  }}
\newcommand{\makelastpage}{\medskip\@contact}

\newcommand{\DUnion}{\bigsqcup\limits}

\newcommand{\Union}{\bigcup\limits}
\newcommand{\inter}{\cap}

\newcommand{\C}{\mathbb{C}}

\newcommand{\R}{\mathbb{R}}

\newcommand{\Z}{\mathbb{Z}}

\newcommand{\op}{\mathrm{op}}
\DeclareMathOperator{\Ob}{Ob}
\DeclareMathOperator{\id}{id}

\newcommand{\BDC}{\mathbf{D}^{\mathrm{b}}}
\newcommand{\PDC}{\mathbf{D}^+}
\newcommand{\NDC}{\mathbf{D}^-}
\newcommand{\TDC}{\mathbf{D}}

\newcommand{\DSum}{\bigoplus}

\newcommand{\ilim}[1][]{\mathop{\varinjlim}\limits_{#1}}

\renewcommand{\to}[1][]{\xrightarrow[#1]{}}
\newcommand{\from}[1][]{\xleftarrow[#1]{}}
\newcommand{\isoto}[1][]{\xrightarrow[#1]{\sim}}

\newcommand{\Endo}[1][]{\mathrm{End}_{\raise1.5ex\hbox to.1em{}#1}}

\newcommand{\Hom}[1][]{\mathrm{Hom}_{\raise1.5ex\hbox to.1em{}#1}}

\newcommand{\RHom}[1][]{\mathrm{RHom}_{\raise1.5ex\hbox to.1em{}#1}}

\newcommand{\Ext}[2][]{\mathrm{Ext}_{\raise1.5ex\hbox to.1em{}#1}^{#2}}

\newcommand{\THom}[1][]{\mathrm{THom}_{\raise1.5ex\hbox to.1em{}#1}}

\newcommand{\Tens}[1][]{\mathbin{\otimes_{\raise1.5ex\hbox to-.1em{}#1}}}

\newcommand{\LTens}[1][]{\mathbin{\otimes_{\raise1.5ex\hbox to-.1em{}#1}^{L}}}

\newcommand{\Tor}[2][]{\mathrm{Tor}^{\raise1.5ex\hbox to.1em{}#1}_{#2}}

\def\sha{\mathcal{A}}
\def\shb{\mathcal{B}}
\def\shc{\mathcal{C}}

\def\she{\mathcal{E}}
\def\shf{\mathcal{F}}
\def\shg{\mathcal{G}}

\def\shi{\mathcal{I}}

\def\shk{\mathcal{K}}
\def\shl{\mathcal{L}}
\def\shm{\mathcal{M}}
\def\shn{\mathcal{N}}

\def\shp{\mathcal{P}}

\def\shr{\mathcal{R}}

\def\shu{\mathcal{U}}
\def\shv{\mathcal{V}}

\newcommand{\sect}{\varGamma}

\renewcommand{\hom}[1][]{{\mathcal{H}om}_{\raise1.5ex\hbox to.1em{}#1}}

\newcommand{\rhom}[1][]{{R\mathcal{H}om}_{\raise1.5ex\hbox to.1em{}#1}}

\newcommand{\ext}[2][]{{\mathcal{E}xt}_{\raise1.5ex\hbox to.1em{}#1}^{#2}}

\newcommand{\thom}[1][]{{T\mathcal{H}om}_{\raise1.5ex\hbox to.1em{}#1}}

\newcommand{\tens}[1][]{\mathbin{\otimes_{\raise1.5ex\hbox to-.1em{}#1}}}

\newcommand{\ltens}[1][]{\mathbin{\otimes_{\raise1.5ex\hbox to-.1em{}#1}^{L}}}

\newcommand{\tor}[2][]{{\mathcal{T}or}^{\raise1.5ex\hbox to.1em{}#1}_{#2}}

\newcommand{\wtens}{\mathbin{\mathop{\otimes}\limits^{{}_{\mathrm{w}}}}}

\DeclareMathOperator{\supp}{supp}

\newcommand{\oim}[1]{{#1}_*}
\newcommand{\eim}[1]{{#1}_!}

\newcommand{\roim}[1]{{R#1}_*}
\newcommand{\reim}[1]{{R#1}_!}

\newcommand{\opb}[1]{#1^{-1}}

\newcommand{\epb}[1]{#1^{!}}

\newcommand{\GHom}[1][]{\mathrm{GHom}_{\raise1.5ex\hbox to.1em{}#1}}

\newcommand{\GExt}[2][]{\mathrm{GExt}_{\raise1.5ex\hbox to.1em{}#1}^{#2}}

\newcommand{\FHom}[1][]{\mathrm{FHom}_{\raise1.5ex\hbox to.1em{}#1}}

\newcommand{\ghom}[1][]{{\mathcal{GH}om}_{\raise1.5ex\hbox to.1em{}#1}}

\newcommand{\gext}[2][]{{\mathcal{GE}xt}_{\raise1.5ex\hbox to.1em{}#1}^{#2}}

\newcommand{\fhom}[1][]{{\mathcal{FH}om}_{\raise1.5ex\hbox to.1em{}#1}}

\newcommand{\tenstop}[1][]{\mathbin{\hat{\otimes}_{\raise1.5ex\hbox to-.1em{}#1}}}

\newcommand{\homtop}[1][]{\mathcal{L}_{\raise1.5ex\hbox to.1em{}#1}}

\newcommand{\Homtop}[1][]{\mathrm{L}_{\raise1.5ex\hbox to.1em{}#1}}

\newcommand{\D}{\mathcal{D}}

\renewcommand{\O}{\mathcal{O}}

\DeclareMathOperator{\chv}{char}

\def\absdoim#1{\underline{#1}_*}
\def\reldoim[#1]#2{\underline{#2}_{|{#1}*}}
\def\doim{\@ifnextchar [{\reldoim}{\absdoim}}

\def\absdeim#1{\underline{#1}_*}
\def\reldeim[#1]#2{\underline{#2}_{|{#1}*}}
\def\deim{\@ifnextchar [{\reldeim}{\absdeim}}

\def\absdopb#1{\underline{#1}^{-1}}
\def\reldopb[#1]#2{\underline{#2}_{|{#1}}^{-1}}
\def\dopb{\@ifnextchar [{\reldopb}{\absdopb}}

\def\absboim#1{\underline{\underline{#1}}_*}
\def\relboim[#1]#2{\underline{\underline{#2}}_{|{#1}*}}
\def\boim{\@ifnextchar [{\relboim}{\absboim}}

\def\absbeim#1{\underline{\underline{#1}}_*}
\def\relbeim[#1]#2{\underline{\underline{#2}}_{|{#1}*}}
\def\beim{\@ifnextchar [{\relbeim}{\absbeim}}

\def\absbopb#1{\underline{\underline{#1}}^*}
\def\relbopb[#1]#2{\underline{\underline{#2}}_{|{#1}}^*}
\def\bopb{\@ifnextchar [{\relbopb}{\absbopb}}

\newtheorem{theorem}{Theorem}[section]

\newtheorem{proposition}[theorem]{Proposition}
\newtheorem{lemma}[theorem]{Lemma}
\newtheorem{conjecture}[theorem]{Conjecture}

\theorembodyfont{\rmfamily}

\newtheorem{definition}[theorem]{Definition}
\newtheorem{example}[theorem]{Example}
\newtheorem{remark}[theorem]{Remark}

\numberwithin{equation}{section}

\renewcommand{\to}[1][]{\xrightarrow[]{#1}}
\renewcommand{\from}[1][]{\xleftarrow[]{#1}}

\newcommand{\approxto}{\xrightarrow{\smash{\approx}}}

\newcommand{\Ring}{A}
\newcommand{\CRing}{R}

\newcommand{\univU}{{\mathcal{U}}}

\newcommand{\Aut}[1][]{\mathrm{Aut}_{\raise1.5ex\hbox to.1em{}#1}}

\renewcommand{\BDC}{\mathsf{D}^{\mathrm{b}}}
\renewcommand{\PDC}{\mathsf{D}^+}
\renewcommand{\NDC}{\mathsf{D}^-}
\renewcommand{\TDC}{\mathsf{D}}

\newcommand{\category}{\mathsf}

\newcommand{\catc}{{\mathsf{C}}}
\newcommand{\catd}{{\mathsf{D}}}

\newcommand{\catu}{{\mathsf{U}}}

\newcommand{\catx}{{\mathsf{X}}}

\newcommand{\catHom}[1][]{\mathsf{Hom}_{\raise1.5ex\hbox to.1em{}#1}}
\newcommand{\catEnd}[1][]{\mathsf{End}_{\raise1.5ex\hbox to.1em{}#1}}
\newcommand{\catMod}{\mathsf{Mod}}

\newcommand{\ring}{\sha}
\newcommand{\ringi}{\shb}
\newcommand{\ringii}{\shc}
\newcommand{\cring}{\shr}

\newcommand{\shHom}[1][]{{\mathcal{H}om}_{\raise1.5ex\hbox 
to.1em{}#1}}
\newcommand{\shEndo}[1][]{{\mathcal{E}nd}_{\raise1.5ex\hbox 
to.1em{}#1}}
\newcommand{\shAut}[1][]{{\mathcal{A}ut}_{\raise1.5ex\hbox 
to.1em{}#1}}
\newcommand{\shInn}[1][]{{\mathcal{I}nn}_{\raise1.5ex\hbox 
to.1em{}#1}}
\newcommand{\shOut}[1][]{{\mathcal{O}ut}_{\raise1.5ex\hbox 
to.1em{}#1}}

\newcommand{\prstkp}{\mathfrak{P}}
\newcommand{\prstkq}{\mathfrak{Q}}

\newcommand{\stkm}{\mathfrak{M}}
\newcommand{\stkn}{\mathfrak{N}}

\newcommand{\stks}{\mathfrak{S}}
\newcommand{\stkt}{\mathfrak{T}}

\newcommand{\stkPSh}{\mathfrak{PSh}}
\newcommand{\stkSh}{\mathfrak{Sh}}
\newcommand{\stkHom}[1][]{\mathfrak{Hom}_{\raise1.5ex\hbox 
to.1em{}#1}}
\newcommand{\stkMod}{\mathfrak{Mod}}

\newcommand{\twid}{\mathsf{1}}
\newcommand{\twst}{\mathsf{t}}
\newcommand{\twsu}{\mathsf{u}}

\newcommand{\cringalg}{\cring\text{-alg}}

\newcommand{\Pic}{\operatorname{Pic}}
\newcommand{\TWS}{\operatorname{Tw}}
\newcommand{\TDO}{\operatorname{TDO}}
\newcommand{\Br}{\operatorname{Br}}

\newcommand{\ad}{{\operatorname{ad}}}
\newcommand{\rightad}{{\operatorname{r}}}

\newcommand{\shift}[1]{[#1]}
\newcommand{\twist}[1]{{\scriptstyle{(#1)}}}

\newcommand{\ctens}{\mathbin{\mathop{\otimes}\limits_{\C}}}
\newcommand{\otens}{\mathbin{\mathop{\otimes}\limits_{\O}}}

\newcommand{\lotens}{\mathbin{\mathop{\otimes}\limits_{\O}^{L}}}

\newcommand{\oopb}[1]{#1^{*}}

\renewcommand{\oim}[1]{#1_*}
\renewcommand{\eim}[1]{#1_!}

\newcommand{\DD}{\mathbb{D}}
\renewcommand{\deim}[1]{{\DD #1}_{!}}
\renewcommand{\doim}[1]{{\DD #1}_{*}}
\renewcommand{\dopb}[1]{\DD #1^{*}}

\newcommand{\dtens}{\mathbin{\mathop{\otimes}\limits^{\DD}}}
\newcommand{\dhom}[1][]{{\DD\mathcal{H}om}_{\raise1.5ex\hbox 
to.1em{}#1}}

\newcommand{\dcirc}{\mathbin{\mathop{\circ}\limits^{\DD}}}

\newcommand{\AX}{\sha}
\newcommand{\BX}{\sha'}
\newcommand{\AY}{\shb}

\newcommand{\AXY}{\tdoinv\AX\tdoetens\AY}

\newcommand{\OAX}[1][]{\O_{\AX#1}}

\newcommand{\OAXY}{\O_{\AXY}}

\newcommand{\CX}{\C_X}
\newcommand{\CY}{\C_Y}

\newcommand{\OX}{\O_X}
\newcommand{\OvX}{\Omega_X}
\newcommand{\OY}{\O_Y}

\newcommand{\Ovf}{\Omega_f}
\newcommand{\DX}{\D_X}

\newcommand{\PP}{\mathbb{P}}
\newcommand{\OP}{\O_\PP}

\newdimen\wdsquare \newdimen\wdsquares \newdimen\wdsquaress
\newdimen\wdsquastkern \newdimen\wdsquastkerns \newdimen\wdsquastkernss 
\newdimen\wdsqushrpkern \newdimen\wdsqushrpkerns \newdimen\wdsqushrpkernss 
\newdimen\htsquare \newdimen\htsquares \newdimen\htsquaress
\newdimen\htsquastkern \newdimen\htsquastkerns \newdimen\htsquastkernss
\newdimen\htsqushrpkern \newdimen\htsqushrpkerns \newdimen\htsqushrpkernss

\setbox0=\hbox{$\square$}
\wdsquare=\wd0
\htsquare=\ht0
\setbox0=\hbox{$\ast$}
\advance\wdsquastkern by -.5\wd0
\advance\htsquastkern by -.5\ht0
\setbox0=\hbox{$\scriptstyle\sharp$}
\advance\wdsqushrpkern by -.5\wd0
\advance\htsqushrpkern by -.5\ht0
\advance\htsqushrpkern by .5\dp0

\setbox0=\hbox{$\scriptstyle\square$}
\wdsquares=\wd0
\htsquares=\ht0
\setbox0=\hbox{$\scriptstyle\ast$}
\advance\wdsquastkerns by -.5\wd0
\advance\htsquastkerns by -.5\ht0
\setbox0=\hbox{$\scriptscriptstyle\sharp$}
\advance\wdsqushrpkerns by -.5\wd0
\advance\htsqushrpkerns by -.5\ht0
\advance\htsqushrpkerns by .5\dp0

\setbox0=\hbox{$\scriptscriptstyle\square$}
\wdsquaress=\wd0
\htsquaress=\ht0
\setbox0=\hbox{$\scriptscriptstyle\ast$}
\advance\wdsquastkernss by -.5\wd0
\advance\htsquastkernss by -.5\ht0
\setbox0=\hbox{$\scriptscriptstyle\sharp$}
\advance\wdsqushrpkernss by -.5\wd0
\advance\htsqushrpkernss by -.5\ht0
\advance\htsqushrpkernss by .5\dp0

\newcommand{\stketimesi}{\hbox to\wdsquare{$\square$\kern\wdsquastkern%
\raise\htsquastkern\hbox{$\ast$}\hfil}}

\newcommand{\stketimesii}{\hbox to\wdsquares{$\scriptstyle\square$\kern\wdsquastkerns%
\raise\htsquastkerns\hbox{$\scriptstyle\ast$}\hfil}}

\newcommand{\stketimesiii}{\hbox to\wdsquaress{$\scriptscriptstyle\square$\kern\wdsquastkernss%
\raise\htsquastkernss\hbox{$\scriptscriptstyle\ast$}\hfil}}

\newcommand{\stketimes}{\mathchoice{\stketimesi}{\stketimesi}{\stketimesii}{\stketimesiii}}

\newcommand{\tdoetimesi}{\hbox to\wdsquare{$\square$\kern\wdsqushrpkern%
\raise\htsqushrpkern\hbox{$\scriptstyle\sharp$}\hfil}}

\newcommand{\tdoetimesii}{\hbox to\wdsquares{$\scriptstyle\square$\kern\wdsqushrpkerns%
\raise\htsqushrpkerns\hbox{$\scriptscriptstyle\sharp$}\hfil}}

\newcommand{\tdoetimesiii}{\hbox to\wdsquaress{$\scriptscriptstyle\square$\kern\wdsqushrpkernss%
\raise\htsqushrpkernss\hbox{$\scriptscriptstyle\sharp$}\hfil}}

\newcommand{\tdoetimes}{\mathchoice{\tdoetimesi}{\tdoetimesi}{\tdoetimesii}{\tdoetimesiii}}

\newcommand{\tdotimes}{\mathchoice%
{\setbox0=\hbox{$\sharp$}\hbox{\raise.5\dp0\hbox{$\sharp$}}}
{\setbox0=\hbox{$\sharp$}\hbox{\raise.5\dp0\hbox{$\sharp$}}}
{\setbox0=\hbox{$\scriptstyle\sharp$}\hbox{\raise.5\dp0\hbox{$\scriptstyle\sharp$}}}
{\setbox0=\hbox{$\scriptscriptstyle\sharp$}\hbox{\raise.5\dp0\hbox{$\scriptscriptstyle\sharp$}}}
}

\newcommand{\stktimes}{\circledast}

\newcommand{\stktens}[1][\cring]{\mathbin{\stktimes_{\raise1.5ex\hbox 
to-.1em{}#1}}}
\newcommand{\stketens}{\mathbin{\stketimes}}
\newcommand{\stkopb}[1]{#1^{\stktimes}}
\newcommand{\stkinv}[2][-1]{#2^{\stktimes#1}}

\newcommand{\tdotens}{\mathbin{\tdotimes}}
\newcommand{\tdoetens}{\mathbin{\tdoetimes}}
\newcommand{\tdoopb}[1]{#1^{\sharp}}
\newcommand{\tdoinv}[2][-1]{#2^{\sharp#1}}

\contact[dagnolo@math.unipd.it, pietro@math.jussieu.fr]{A.D'A.: Dipartimento di Matematica Pura ed
Applicata; Universit\`a di Padova; Via G.~Belzoni, 7; 35131 Padova;
Italy\\
P.P.: Dipartimento di Matematica Pura ed
Applicata; Universit\`a di Padova; Via G.~Belzoni, 7; 35131 Padova;
Italy; and: Analyse Alg\'ebrique; Institut de Math\'ematiques; 175,
rue du Chevaleret; 75013 Paris; France}

\begin{document}

\author{Andrea D'Agnolo\thanks{A.D'A.\ had the occasion to visit the
Research Institute for Mathematical Sciences of Kyoto University
during the preparation of this paper.  Their warm hospitality is
gratefully acknowledged.} \and Pietro Polesello\thanks{P.P.\ was
partially supported by INdAM during the preparation of this paper}}

\title{Stacks of twisted modules and integral transforms}

\maketitle

\begin{abstract}
Stacks were introduced by Grothendieck and Giraud and are, roughly 
speaking, sheaves of categories.  Kashiwara developed the theory of 
twisted modules, which are objects of stacks locally equivalent to 
stacks of modules over
sheaves of rings.  In this paper we recall these notions, and we develop the 
formalism of operations for stacks of twisted modules.
As an application, we state a twisted version of an adjunction 
formula which is of use in
the theory of integral transforms for sheaves and $D$-modules.
\end{abstract}

\begin{classification} 14A20, 32C38, 35A22
\end{classification}

\tableofcontents

\section*{Introduction}

Stacks are, roughly speaking, sheaves of categories.  They were
introduced by Grothendieck and Giraud~\cite{Giraud1971} 
in algebraic geometry where
some special stacks, called gerbes, are now commonly used in moduli
problems to describe objects with automorphisms (see for
example~\cite{Behrend2002,Laumon-Moret-Bailly2000}).  Recently, gerbes
have infiltrated differential geometry and mathematical physics (see
for 
example~\cite{Brylinski1993,Murray1996,Hitchin1999,Breen-Messing2001}).

We are interested here in the related notion of twisted modules, which
are objects of stacks locally equivalent to stacks of modules over
sheaves of rings. The simplest example is that of stacks of twisted $\cring$-modules on a locally ringed space $(X,\cring)$. These can be considered as higher cohomological analogues to line bundles. More precisely, line bundles on $X$ are sheaves of $\cring$-modules locally isomorphic to $\cring$, and their isomorphism classes describe the cohomology group $H^1(X;\cring^\times)$. 
Stacks of twisted $\cring$-modules are $\cring$-linear stacks on $X$ locally 
equivalent to the stack $\stkMod(\cring)$ of $\cring$-modules and, as we shall recall, their equivalence classes describe the cohomology group $H^2(X;\cring^\times)$. As line bundles correspond to principal $\cring^\times$-bundles, so stacks of twisted $\cring$-modules correspond to gerbes with band $\cring^\times$. However, this correspondence no longer holds for the more general type of stacks of twisted modules that we  consider here.

Twisted modules appear in works by Kashiwara on
representation theory~\cite{Kashiwara1989} and on
quantization~\cite{Kashiwara1996}.  In the first case, they were used to
describe solutions on flag manifolds to quasi-equivariant modules
over rings of twisted differential operators (see 
also~\cite{Kashiwara-Schmid1994}).  In the second
case, twisted modules turned out to be the natural framework for a global study of
microdifferential systems on a holomorphic contact manifold (see 
also~\cite{Kontsevich2001,Polesello-Schapira}). 
Rings of microdifferential operators can be locally defined on a contact manifold, 
but do not necessarily exist globally. Kashiwara proved that there 
exists a globally defined $\C$-linear stack which is locally equivalent to 
the stack of modules over a ring of microdifferential operators.

Twisted modules induced by Azumaya algebras are used in~\cite{Caldararu2002,Donagi-Pantev2003} in relation with the Fourier-Mukai transform.

Motivated by Kashiwara's work on quantization, we consider here twisted modules over sheaves of rings which are not necessarily commutative nor globally defined. More precisely, let $X$ be a topological space, or more generally a site, and $\cring$ a sheaf of commutative rings on $X$. Then $\stkm$ is a stack of $\cring$-twisted modules on $X$ if it is $\cring$-linear and there exist an open covering $X = \Union\nolimits_{i\in I} U_i$, sheaves of $\cring|_{U_i}$-algebras $\ring_i$, and $\cring|_{U_i}$-equivalences $\stkm|_{U_i} \approx \stkMod(\ring_i)$, where $\stkMod(\ring_i)$ denotes the stack of left $\ring_i$-modules on $U_i$.  We review the notions of stack and stack of twisted modules in Section~\ref{se:stacks}, restricting to the case of topological spaces for simplicity of exposition. 

Morita theory is the basic tool to deal with stacks of $\cring$-twisted modules, and we use it to develop the formalism of operations, namely duality $\stkinv{(\cdot)}$, internal product $\stktens$, and inverse image $\stkopb f$ by a continuous map $f\colon Y \to X$. If $\ring$ and $\ring'$ are sheaves of $\cring$-algebras on $X$, these operations satisfy $\stkinv{\stkMod(\ring)} \approx \stkMod(\ring^\op)$, $\stkMod(\ring) \stktens \stkMod(\ring') \approx \stkMod(\ring \tens[\cring] \ring')$, and $\stkopb f \stkMod(\ring) \approx \stkMod(\opb f\ring)$. With this formalism at hand, we then describe Grothendieck's six operations for derived categories of twisted modules over locally compact Hausdorff topological spaces. This is the content of Section~\ref{se:operations}.

Morita theory is used again in Section~\ref{se:descent} to describe 
effective descent data attached to semisimplicial complexes.  In particular, we get a Cech-like classification of stacks of $\cring$-twisted modules, 
with invertible bimodules as cocycles, which is parallel to the bitorsor description of gerbes in~\cite{Breen-Messing2001}.

In Section~\ref{se:local}, assuming that $\cring$ is a commutative local 
ring, we recall the above mentioned classification of stacks of twisted 
$\cring$-modules in terms of $H^2(X;\cring^\times)$. We then consider the case of twisted modules induced by ordinary modules over an inner form of a given 
$\cring$-algebra.  This allows us to present in a unified manner the 
examples provided by modules over Azumaya algebras and over rings of
twisted differential operators.  Finally,
we state a twisted version of an adjunction formula for sheaves and
$D$-modules, which is of use in the theory of
integral transforms with regular kernel, as the Radon-Penrose
transform.

This paper is in part a survey and in part original. The survey 
covers material from Kashiwara's papers~\cite{Kashiwara1989,Kashiwara1996}, from his joint works~\cite{Kashiwara-Schmid1994,Kashiwara-Tanisaki1996}, and from the last chapter of his forthcoming book with Pierre Schapira~\cite{Kashiwara-Schapira}. 
The main original contribution is in establishing the formalism of operations for stacks of twisted modules.

It is a pleasure to thank Masaki Kashiwara for several useful 
discussions and insights. We also wish to thank him and Pierre 
Schapira for allowing us to use results from a preliminary version 
of their book~\cite{Kashiwara-Schapira}.

\section{Stacks of twisted modules}\label{se:stacks}

The theory of stacks is due to Grothendieck and Giraud~\cite{Giraud1971}.  We review it here restricting for simplicity to the case of stacks on topological spaces (thus avoiding the notions of site 
and of fibered category).  Finally, we recall the notion of stack of twisted modules, considering the case of modules over rings which are not necessarily commutative nor globally defined.
Our main references were~\cite{Kashiwara1989,Kashiwara1996,Kashiwara-Schapira2001,Kashiwara-Schapira}.

\subsection{Prestacks}

We assume that the reader is familiar with the basic notions of
category theory, as those of category, functor between categories,
transformation between functors (also called morphism of functors),
and equivalence of categories.

If $\catc$ is a category, we denote by $\Ob(\catc)$ the
set\footnote{Following Bourbaki's appendix in~\cite{SGA4}, one way to 
avoid the paradoxical
situation of dealing with the set of all sets is to consider
universes, which are ``big'' sets of sets stable by most of the
set-theoretical operations.  We assume here to be given a universe
$\univU$ and, unless otherwise stated, all categories $\catc$ will be
assumed to be $\univU$-categories, i.e.\ categories such that
$\Ob(\catc)\subset\univU$ and $\Hom[\catc](c,d) \in \univU$ for every
pair of objects.} of its objects, and by $\Hom[\catc](c,d)$ the set of
morphisms between the objects $c$ and $d$.  The identity of 
$\Hom[\catc](c,c)$ will be denoted by $\id_c$.

Denote by $\catc^\op$ the opposite category, which has the same
objects as $\catc$ and reversed morphisms $\Hom[\catc^\op](c,d) =
\Hom[\catc](d,c)$.  If $\catd$ is another category, denote by
$\catHom(\catc,\catd)$ the category of functors from
$\catc$ to $\catd$, with transformations as morphisms.

Let $X$ be a topological space, and denote by $\catx$ the category of
its open subsets with inclusion morphisms.  Recall that the category
of presheaves on $X$ with values in a category $\catc$ is the category
$\catHom(\catx^\op,\catc)$ of contravariant functors from $\catx$ to
$\catc$. In particular, presheaves of sets are obtained by taking $\catc = \category{Set}$,
the category of sets\footnote{More
precisely, $\category{Set}$ denotes the $\univU$-category of sets
belonging to the fixed universe $\univU$.} and maps of sets. 

Considering $\catc = \category{Cat}$, the category of categories\footnote{More precisely, let
$\shv$ be another fixed universe with $\shu\in\shv$.  Then
$\category{Cat}$ denotes the $\shv$-category whose objects are
$\shu$-categories.  From now on we will leave to the reader who feels
that need the task of making the universes explicit.} and functors, one has a notion of presheaf of
categories. This a functor $\shf\colon \catx^\op \to
\category{Cat}$, and if $W\to[v]V\to[u]U$ are inclusions
of open sets, the restriction functors $\shf(u)\colon \shf(U) \to
\shf(V)$ and $\shf(v)\colon \shf(V) \to \shf(W)$ are thus required to
satisfy the equality $\shf(v) \circ \shf(u) = \shf(u \circ v)$.  Such a
requirement is often too strong in practice, and the notion of
prestack is obtained by weakening this equality to an isomorphism of
functors, i.e.\ to an invertible transformation.  

In other words, prestacks are the 2-categorical\footnote{Roughly speaking, a 2-category (refer to
\cite[\S9]{Street1996} for details) $\mathbf{C}$ is a ``category
enriched in $\category{Cat}$'', i.e.\ a category whose morphism sets
are the object sets of categories $\catHom[\mathbf{C}](c,d)$, such
that composition is a functor.  Morphisms in the category
$\catHom[\mathbf{C}](c,d)$ are called 2-cells.  The basic example is
the 2-category $\mathbf{Cat}$ which has categories as objects,
functors as morphisms, and transformations as 2-cells.

There is a natural notion of pseudo-functor between 2-categories,
preserving associativity for the composition functor only up to an
invertible 2-cell.  Then a prestack (see~\cite[expos\'e VI]{SGA1}) is 
a pseudo-functor
$\mathbf{X}^\op \to \mathbf{Cat}$, where $\mathbf{X}^\op$ is the
2-category obtained by trivially enriching $\catx^\op$ with identity
2-cells.  Functors of prestacks and their transformations are
transformations and modifications of pseudo-functors, respectively.

Note that Corollary 9.2 of~\cite{Street1996} asserts that any prestack
is equivalent, in the 2-category of pseudo-functors, to a presheaf of
categories.  However, this equivalence is not of practical use for our
purposes.} version of
presheaves of categories.  However, we prefer not to use the language of
$2$-categories, giving instead the unfolded definition of prestack.

\begin{definition}
A prestack $\prstkp$ on $X$ consists of the data
\begin{itemize}
    \item[(a)] for every open subset $U\subset X$, a category
    $\prstkp(U)$,

    \item[(b)] for every inclusion $V \to[u] U$ of open subsets, a
    functor $\prstkp(u)\colon \prstkp(U)\to\prstkp(V)$, called
    restriction functor, and

    \item[(c)] for every inclusion $W \to[v] V \to[u] U$ of open
    subsets, invertible transformations $\prstkp(v,u)\colon
    \prstkp(v)\circ\prstkp(u)\Rightarrow\prstkp(u\circ v)$ of 
functors from
    $\prstkp(U)$ to $\prstkp(W)$,
\end{itemize}
subject to the conditions
\begin{itemize}
    \item[(i)] if $U$ is an open subset, then
    $\prstkp(\id_U)=\id_{\prstkp(U)}$ and $\prstkp(\id_U, \id_U) =
    \id_{\id_{\prstkp(U)}}$;
    
    \item[(ii)] if $Y \to[w] W \to[v] V \to[u] U$ are inclusions of
    open subsets, then the following diagram of functors from
    $\prstkp(U)$ to $\prstkp(Y)$ commutes
$$
\xymatrix@C=6em{ \prstkp(w) \circ \prstkp(v) \circ \prstkp(u)
\ar@{=>}[r]^-{\prstkp(w,v) \circ \id_{\prstkp(u)}} 
\ar@{=>}[d]^{\id_{\prstkp(w)} \circ \prstkp(v,u)} &
\prstkp(v\circ w) \circ \prstkp(u) \ar@{=>}[d]^{\prstkp(v\circ w,u)}
\\
\prstkp(w) \circ \prstkp(u\circ v) \ar@{=>}[r]^-{\prstkp(w,u\circ v)}
& \prstkp(u\circ v \circ w).  }
$$
In particular, $\prstkp(u,\id_U) = \prstkp(\id_{V},u) = 
\id_{\prstkp(u)}$.
\end{itemize}
\end{definition}

For $\shf\in \prstkp(U)$ and $V\to[u] U$ an open inclusion, one
usually writes $\shf\vert_{V}$ instead of $\prstkp(u)(\shf)$.  One
denotes by $\prstkp\vert_U$ the natural restriction of $\prstkp$ to
$U$ given by $V\mapsto \prstkp(V)$ for $V\subset U$.

\begin{definition}
Let $\prstkp$ and $\prstkq$ be prestacks on $X$.  A functor of
prestacks $\varphi \colon \prstkp \to \prstkq$ consists of the data
\begin{itemize}
    \item[(a)] for any open subset $U\subset X$, a functor $\varphi(U)
    \colon \prstkp(U) \to \prstkq(U)$,

    \item[(b)] for any open inclusion $V \to[u] U$, an invertible
    transformation $\varphi(u) \colon \varphi(V) \circ \prstkp(u)
    \Rightarrow \prstkq(u) \circ \varphi(U)$ of functors from
    $\prstkp(U)$ to $\prstkq(V)$,
\end{itemize}
subject to the condition
\begin{itemize}
     \item[(i)] if $W \to[v] V \to[u] U$ are inclusions of open
     subsets, then the following diagram of functors from $\prstkp(U)$
     to $\prstkq(W)$ commutes
$$
\!\!\!\!\!\!\!
\xymatrix@C=3.5em{ \varphi(W) \circ \prstkp(v) \circ \prstkp(u)
\ar@{=>}[r]^-{\varphi(v) \circ \id_{\prstkp(u)}} 
\ar@{=>}[d]^{\id_{\varphi(W)} \circ \prstkp(v,u)} & \prstkq(v)
\circ \varphi(V) \circ \prstkp(u) \ar@{=>}[r]^-{\id_{\prstkq(v)} 
\circ \varphi(u)} &
\prstkq(v) \circ \prstkq(u) \circ \varphi(U)
\ar@{=>}[d]_{\prstkq(v,u) \circ \id_{\varphi(U)}} \\
\varphi(W) \circ \prstkp(u\circ v) \ar@{=>}[rr]^-{\varphi(u\circ v)} &
& \prstkq(u\circ v) \circ \varphi(U) .  } $$
In particular, $\varphi(\id_U) = \id_{\varphi(U)}$.
\end{itemize}
\end{definition}

\begin{definition}
Let $\varphi, \psi \colon \prstkp \to \prstkq$ be functors of
prestacks.  A transformation $\alpha\colon \varphi \Rightarrow \psi$
of functors of prestacks consists of the data
\begin{itemize}
    \item[(a)] for any open subset $U\subset X$, a transformation
    $\alpha(U) \colon \varphi(U) \Rightarrow \psi(U)$ of functors from
    $\prstkp(U)$ to $\prstkq(U)$,
\end{itemize}
such that
\begin{itemize}
     \item[(i)] if $V\to[u] U$ is an inclusion of open subsets, then
     the following diagram of functors from $\prstkp(U)$ to
     $\prstkq(V)$ commutes
$$
\xymatrix@C=4em{ \varphi(V) \circ \prstkp(u) \ar@{=>}[r]^-{\alpha(V) 
\circ \id_{\prstkp(u)}}
\ar@{=>}[d]^{\varphi(u)} & \psi(V) \circ \prstkp(u)
\ar@{=>}[d]^{\psi(u)} \\
\prstkq(u) \circ \varphi(U) \ar@{=>}[r]^-{\id_{\prstkq(u)} \circ 
\alpha(U)} & \prstkq(u)
\circ \psi(U).  }
$$
\end{itemize}
\end{definition}

An example of prestack is the prestack $\stkPSh_X$ of presheaves of sets on $X$. It associates to an open subset $U\subset X$ the category 
$\catHom(\catu^\op,\category{Set})$ of presheaves of sets on $U$, and to an open inclusion $V 
\subset U$ the restriction functor $\stkPSh_X(U) \to \stkPSh_X(V)$, 
$\shf \mapsto \shf\vert_V$. For open inclusions $W \subset V \subset 
U$ one has $\shf\vert_V\vert_W = \shf\vert_W$, so that $\stkPSh_X$ is 
in fact a presheaf of categories.

If $\prstkp$ and $\prstkq$ are prestacks, one gets another prestack 
$\stkHom(\prstkp, \prstkq)$ by associating to an open subset 
$U\subset X$ the category $\catHom(\prstkp\vert_U, \prstkq\vert_U)$ 
of functors of prestacks from
$\prstkp\vert_U$ to $\prstkq\vert_U$, with transformations of functors
of prestacks as morphisms, and with the natural restriction functors. Note that $\stkHom(\prstkp, \prstkq)$ is actually a presheaf of categories. 

(Pre)stacks which are not (pre)sheaves of categories will 
appear in Section~\ref{se:operations}.

\subsection{Stacks}

The analogy between presheaves and prestacks goes on for sheaves and
stacks.

Let $X$ be a topological space.  Given a family of subsets
$\{U_i\}_{i\in I}$ of $X$, let us use the notations
$$
U_{ij} = U_{i} \cap U_{j}, \quad U_{ijk} = U_{i} \cap U_{j} \cap
U_{k}, \quad \text{etc.}
$$
Recall that a presheaf of sets $\shf$ on $X$ is called a sheaf if for
any open subset $U\subset X$, and any open covering $\{U_i\}_{i\in I}$
of $U$, the natural sequence given by the restriction maps
$$
\xymatrix{ \shf(U) \ar[r] & \prod_{i\in I} \shf(U_i) \ar@<.2pc>[r]
\ar@<-.2pc>[r] & \prod_{i,j\in I} \shf(U_{ij}) }
$$
is exact, i.e.\ if for any family of sections $s_i\in\shf(U_i)$
satisfying $s_i|_{U_{ij}} = s_j|_{U_{ij}}$ there is a unique section
$s\in\shf(U)$ such that $s|_{U_i} = s_i$.

Similarly to the definition of sheaf, a prestack $\stks$ on $X$ is 
called a stack if for any open
subset $U\subset X$, and any open covering $\{U_i\}_{i\in I}$ of $U$,
the natural sequence given by the restriction functors
$$
\xymatrix{ \stks(U) \ar[r] & \prod_{i\in I} \stks(U_i) \ar@<.2pc>[r]
\ar@<-.2pc>[r] & \prod_{i,j\in I} \stks(U_{ij}) \ar[r] \ar@<.3pc>[r]
\ar@<-.3pc>[r] & \prod_{i,j,k\in I} \stks(U_{ijk}) }
$$
is exact in the sense of~\cite[expos\'e XIII]{SGA1}, i.e.\ if the 
category $\stks(U)$ is equivalent to the
category whose objects are families of objects $\shf_i$ of
$\stks(U_i)$ and of isomorphisms $\theta_{ij}\colon
\shf_j\vert_{U_{ij}}\isoto \shf_i\vert_{U_{ij}}$ which are compatible
in the triple intersections, in a natural sense.

More explicitly, recall that a descent datum for $\stks$ on $U$ is a
triplet
\begin{equation}
\label{eq:descent}
F= (\{U_{i}\}_{i\in I}, \{\shf_{i}\}_{i\in I}, \{\theta_{ij}\}_{i,j\in
I}),
\end{equation}
where $\{U_{i}\}_{i\in I}$ is an open covering of $U$, $\shf_i \in
\stks(U_i)$, and $\theta_{ij}\colon \shf_j\vert_{U_{ij}}\isoto
\shf_i\vert_{U_{ij}}$ are isomorphisms such that the following diagram
of isomorphisms commutes
$$
\xymatrix{ \shf_{j}\vert_{U_{ijk}} &
\shf_{j}\vert_{U_{ij}}\vert_{U_{ijk}} \ar[l]_{\stks}
\ar[rr]^{\theta_{ij}\vert_{U_{ijk}}} && 
\shf_{i}\vert_{U_{ij}}\vert_{U_{ijk}}
\ar[r]^{\stks} & \shf_{i}\vert_{U_{ijk}} \\
\shf_{j}\vert_{U_{jk}}\vert_{U_{ijk}} \ar[u]^{\stks} &&&&
\shf_{i}\vert_{U_{ik}}\vert_{U_{ijk}} \ar[u]_{\stks} \\
& \shf_{k}\vert_{U_{jk}}\vert_{U_{ijk}} \ar[r]^-{\stks}
\ar[ul]^{\theta_{jk}\vert_{U_{ijk}}} & \shf_{k}\vert_{U_{ijk}} &
\shf_{k}\vert_{U_{ik}}\vert_{U_{ijk}}. 
\ar[ur]_{\theta_{ik}\vert_{U_{ijk}}}
\ar[l]_-{\stks} & } $$
The descent datum $F$ is called effective if there exist $\shf \in
\stks(U)$ and isomorphisms $\theta_i\colon \shf\vert_{U_i}\isoto
\shf_i$ for each $i$, such that the following diagram of isomorphisms
commutes
$$
\xymatrix{ \shf\vert_{U_{j}}\vert_{U_{ij}} 
\ar[d]_{\theta_{j}\vert_{U_{ij}}}
\ar[r]^{\stks} & \shf\vert_{U_{ij}} & \shf\vert_{U_{i}}\vert_{U_{ij}}
\ar[d]_{\theta_{i}\vert_{U_{ij}}} \ar[l]_{\stks} \\ 
\shf_{j}\vert_{U_{ij}}
\ar[rr]^{\theta_{ij}} && \shf_{i}\vert_{U_{ij}}. }
$$

To $\stks$ a prestack on $X$ is attached a bifunctor of prestacks
$$
\shHom[\stks] \colon \stks^\op \times \stks \to \stkPSh_X,
$$
associating to $\shf,\shg\in\stks(U)$ the presheaf of sets on $U\subset X$ given by
$$
\shHom[\stks|_U](\shf,\shg) \colon V \mapsto
\Hom[\stks(V)](\shf\vert_{V},\shg\vert_{V}).
$$ 

\begin{definition}
\label{de:stack}
\begin{itemize}
\item[(i)]
A prestack $\stks$ on $X$ is called separated if for any open subset 
$U$, and any $\shf,\shg\in\stks(U)$, the presheaf 
$\shHom[\stks|_U](\shf,\shg)$ is a sheaf.
\item[(ii)]
A stack is a separated prestack such that any descent datum is
effective.
\item[(iii)]
Functors and transformations of stacks are functors and
transformations of the underlying prestacks, respectively.
\end{itemize}
\end{definition}

For example, the prestack $\stkSh_X$ of sheaves of sets, associating to $U\subset X$ the category of sheaves of sets on $U$, is actually a 
stack.  As another example, if $\stks$ and $\stkt$ are stacks, then 
the prestack $\stkHom(\stks,\stkt)$ is a stack. (Note that both 
$\stkSh_X$ and $\stkHom(\stks,\stkt)$ are in fact sheaves of 
categories.)

\medskip

One says that a functor of stacks $\varphi\colon\stks\to \stkt$ is an
equivalence if there exists a functor $\psi\colon\stkt\to \stks$,
called a quasi-inverse to $\varphi$, and invertible transformations
$\varphi\circ \psi\Rightarrow \id_{\stkt}$ and $\psi\circ
\varphi\Rightarrow \id_{\stks}$.  One says that $\varphi$ admits a
right adjoint if there exists a functor of stacks $\psi\colon\stkt\to
\stks$, called a right adjoint to $\varphi$, and an invertible
transformation $\shHom[\stkt](\varphi(\cdot),\cdot) \Rightarrow
\shHom[\stks](\cdot,\psi(\cdot))$.  Similarly for left adjoint. 
Finally, if $\stkt = \stkSh_X$, one says that $\varphi\colon \stks \to
\stkSh_X$ is representable if there exists $\shf\in\stks(X)$, called a
representative of $\varphi$, and an invertible transformation $\varphi
\Rightarrow \shHom[\stks](\shf,\cdot)$.

\begin{lemma}
For a functor of stacks to be an equivalence (resp.\ to admit a right
or left adjoint, resp.\ to be representable) is a local property.
\end{lemma}

\begin{proof}
Right or left adjoints and representatives are unique up to unique
isomorphisms, and hence glue together globally.  As for equivalences,
assume that $\varphi$ is locally an equivalence.  Then we have to show
that for each open subset $U\subset X$ the functors $\varphi(U)$ are
fully faithful and essentially surjective.  Being fully faithful is a
local property already for separated prestacks.  Assume that
$\varphi(U_i)$ are essentially surjective for a covering
$U=\Union\nolimits_i U_i$.  Let $\shg\in\Ob(\stkt(U))$, and choose
$\shf_i\in\Ob(\stks(U_i))$ with isomorphisms
$\varphi(U_i)(\shf_i)\isoto \shg|_{U_i}$.  Since $\varphi$
is fully faithful, the restriction morphisms of $\shg|_{U_i}$ give descent data for $\shf_i$.  Finally, since $\stks$ is a stack, one gets
$\shf\in\Ob(\stks(U))$ with $\varphi(U)(\shf)\isoto \shg$.
\end{proof}

\subsection{Constructions of stacks}

The forgetful functor, associating to a sheaf of sets its underlying 
presheaf, has a left adjoint, associating a sheaf $\shp^+$ to a 
presheaf $\shp$. There is a similar
construction associating a stack
$\prstkp^+$ to a prestack $\prstkp$. This is done in two steps as 
follows. 
Consider first the separated prestack $\prstkp^a$, with the same
objects as $\prstkp$ and morphisms
$$
\Hom[\prstkp^a(U)](\shf,\shg) = \sect(U;
\shHom[\prstkp|_U]^+(\shf,\shg)).
$$
Then let $\prstkp^+(U)$ be the category whose objects are descent data
for $\prstkp^a$ on $U$, and whose morphisms $(\{U_{i}\}, \{\shf_{i}\},
\{\theta_{ii'}\}) \to (\{V_j\}, \{\shg_j\}, \{\varpi_{jj'}\})$ consist
of morphisms $\varphi_{ji}\colon \shf_i|_{U_i\cap V_j} \to
\shg_j|_{U_i\cap V_j}$ such that $\varpi_{j'j} \circ \varphi_{ji} =
\varphi_{j'i'} \circ \theta_{i'i}$ on $U_{ii'}\cap V_{jj'}$.

Since sheaves of sets form a stack, descent data for sheaves are
effective.  Similarly, it is possible to patch stacks together.  More
precisely, a descent datum for stacks on $X$ is a quadruplet
\begin{equation}
\label{eq:descentstack}
S = (\{U_i\}_{i\in I}, \{\stks_i\}_{i\in I}, \{\varphi_{ij}\}_{i,j\in
I}, \{\alpha_{ijk}\}_{i,j,k\in I}),
\end{equation}
where $\{U_i\}_{i\in I}$ is an open covering of $X$, $\stks_i$ are
stacks on $U_i$, $\varphi_{ij}\colon \stks_j\vert_{U_{ij}}\approxto
\stks_i\vert_{U_{ij}}$ are equivalences of stacks, and
$\alpha_{ijk}\colon \varphi_{ij}\circ \varphi_{jk} \Rightarrow
\varphi_{ik}$ are invertible transformations of functors from
$\stks_k\vert_{U_{ijk}}$ to $\stks_i\vert_{U_{ijk}}$, such that for
any $i,j,k,l \in I$, the following diagram of transformations of
functors from $\stks_l|_{U_{ijkl}}$ to $\stks_i|_{U_{ijkl}}$ commutes
\begin{equation}
\label{eq:phiijcommut}
\xymatrix@C=4em{ {\varphi_{ij}\circ\varphi_{jk}\circ\varphi_{kl}}
\ar@{=>}[r]^-{\alpha_{ijk} \circ \id_{\varphi_{kl}}} 
\ar@{=>}[d]^{\id_{\varphi_{ij}} \circ \alpha_{jkl}}
&{\varphi_{ik}\circ\varphi_{kl}}\ar@{=>}[d]^{\alpha_{ikl}}\\
{\varphi_{ij}\circ\varphi_{jl}}\ar@{=>}[r]^{\alpha_{ijl}}&\varphi_{il}.
}
\end{equation}

\begin{proposition}
\label{pr:patch}
Descent data for stacks are effective, meaning that given a descent
datum for stacks $S$ as in \eqref{eq:descentstack}, there exist a
stack $\stks$ on $X$, equivalences of stacks $\varphi_i\colon
\stks\vert_{U_i}\approxto\stks_i$, and invertible transformations of 
functors $\alpha_{ij} \colon \varphi_{ij}\circ
\varphi_j\vert_{U_{ij}} \Rightarrow \varphi_i\vert_{U_{ij}}$  such 
that $\alpha_{ij}|_{U_{ijk}} \circ \alpha_{jk}|_{U_{ijk}} = 
\alpha_{ik}|_{U_{ijk}} \circ \alpha_{ijk}$.  The stack
$\stks$ is unique up to equivalence.
\end{proposition}

\begin{proof}[Sketch of proof]
For $U\subset X$ open, denote by $\stks(U)$ the category whose objects
are triplets
$$
F= (\{V_{i}\}_{i\in I}, \{\shf_{i}\}_{i\in I}, \{\xi_{ij}\}_{i,j\in
I}),
$$
where $V_i = U\cap U_i$, $\shf_i\in\Ob(\stks_i(V_i))$, and
$\xi_{ij}\colon \varphi_{ij}(\shf_{j}\vert_{V_{ij}}) \to \shf_{i}
\vert_{V_{ij}}$ are isomorphisms such that for $i,j,k\in I$ the
following diagram commutes
$$
\xymatrix{ \varphi_{ij}( \varphi_{jk}( \shf_{k}\vert_{V_{ijk}}))
\ar[r]^-{\alpha_{ijk}} \ar[d]^{\varphi_{ij}(\xi_{jk}\vert_{V_{ijk}})} &
\varphi_{ik}(\shf_{k}\vert_{V_{ijk}}) 
\ar[d]^{\xi_{ik}\vert_{V_{ijk}}} \\
\varphi_{ij}(\shf_{j}\vert_{V_{ijk}}) 
\ar[r]^-{\xi_{ij}\vert_{V_{ijk}}} &
\shf_{i}\vert_{V_{ijk}} .  }
$$
For $G=(\{V_{i}\}, \{\shg_{i}\}, \{\eta_{ij}\})$, a morphisms
$\gamma\colon F \to G$ in $\stks(U)$ consists of morphisms
$\gamma_i\colon \shf_i \to \shg_i$ in $\stks_i(V_i)$ such that the
following diagram commutes
$$
\xymatrix{ \varphi_{ij}( \shf_j\vert_{V_{ij}}) \ar[r]^-{\xi_{ij}}
\ar[d]^{\varphi_{ij}(\gamma_{j}\vert_{V_{ij}})} & \shf_{i}\vert_{V_{ij}} 
\ar[d]^{\gamma_{i}\vert_{V_{ij}}} \\
\varphi_{ij}(\shg_{j}\vert_{V_{ij}}) \ar[r]^-{\eta_{ij}} &
\shg_{i}\vert_{V_{ij}} .  }
$$
Then one checks that the prestack $\stks\colon U \to \stks(U)$ is a
stack satisfying the requirements in the statement.
\end{proof}

\subsection{Operations}

Let us recall the stack-theoretical analogue of internal and external 
operations for sheaves.

Given two stacks $\stks$ and $\stks'$ on $X$, denote by
$\stks\times\stks'$ the prestack $\stks\times\stks'(U)=
\stks(U)\times\stks'(U)$.  This is actually a stack.
We already noticed that the prestack $\stkHom(\stks,\stks')$ is a 
stack.
If $\stks''$ is another stack, there is a natural equivalence
$$
\stkHom(\stks\times\stks',\stks'') \approxto \stkHom(\stks,
\stkHom(\stks',\stks'')).
$$

Let $f\colon Y\to X$ be a continuous map of topological spaces.
If $\stkt$ is a stack on $Y$, denote by $\oim f\stkt$ the prestack 
$\oim
f\stkt(U) = \stkt(\opb f U)$, which is actually a stack.
If $\stks$ is a stack on $X$, denote by $\opb f\stks = 
(f^\sim\stks)^+$ the stack associated with the prestack $f^\sim\stks$ 
defined as follows.
For $V\subset Y$, $f^\sim\stks(V)$ is the category whose objects are 
the disjoint union
$\DUnion\nolimits_{U\colon \opb f U\supset V}\Ob(\stks(U))$, and whose
morphisms are given by
$$
\Hom[f^{\sim}\stks(V)](\shf^{U},\shf^{U'}) = \ilim[U''\colon U''
\subset U \cap U' ,\ \opb f U'' \supset V]
\Hom[\stks(U'')](\shf^{U}|_{U''},\shf^{U'}|_{U''}),
$$
for $\shf^U\in\Ob(\stks(U))$ and $\shf^{U'}\in\Ob(\stks(U'))$.  
There is a natural equivalence
$$
\oim f \stkHom(\opb f\stks, \stkt) \approxto \stkHom(\stks, \oim f 
\stkt).
$$
\subsection{Linear stacks}

As a matter of conventions, in this paper rings are unitary, and ring 
homomorphisms preserve the unit.
If $\CRing$ is a commutative ring, we call $\CRing$-algebra a not 
necessarily commutative ring $\Ring$ endowed with a ring homomorphism 
$\CRing\to\Ring$ whose image is in the center of $\Ring$.

Let $\CRing$ be a commutative ring. An $\CRing$-linear category, that we will call $\CRing$-category for short, is a category $\catc$ whose morphism 
sets are endowed with a structure of $\CRing$-module such that 
composition is $\CRing$-bilinear.  An
$\CRing$-functor is a functor which is $\CRing$-linear at the level of
morphisms.  Transformations of $\CRing$-functors are simply
transformations of the underlying functors.  Note that if $\catd$ is
another $\CRing$-category, the category $\catHom[\CRing](\catc,\catd)$
of $\CRing$-functors and transformations is again an
$\CRing$-category, the $\CRing$-module structure on the sets
of transformations being defined object-wise.

For each $c\in\Ob(\catc)$ the set of endomorphisms $\Endo[\catc](c)$
has a natural structure of $\CRing$-algebra, with product given by 
composition.  In particular, note
that $\CRing$-algebras are identified with $\CRing$-categories with a
single object. Let us denote for short by $\Endo(\id_\catc)$ the 
$\CRing$-algebra $\Endo[{\catEnd[\CRing](\catc)}](\id_\catc)$. It is a
commutative\footnote{Let
$\alpha,\beta\colon\id_\catc\Rightarrow\id_\catc$ be transformations,
and $c\in\Ob(\catc)$.  By definition of transformation, applying
$\alpha$ to the morphism $\beta(c)$ we get a commutative diagram
$$\xymatrix{c\ar[r]_{\alpha(c)} \ar[d]_{\beta(c)} & c
\ar[d]^{\beta(c)} \\ c\ar[r]^{\alpha(c)} & c.}$$  Note that the
natural morphism $\Endo(\id_\catc)\to\Endo[\catc](c)$,
$\alpha\mapsto\alpha(c)$, identifies $\Endo(\id_\catc)$ canonically
with the center of $\Endo(c)$ for each $c\in \Ob(\catc)$} 
$\CRing$-algebra,
called the center of $\catc$.  Note that $\catc$ is an
$\CRing$-category if and only if $\catc$ is a $\Z$-category (also
called preadditive category) endowed with a ring homomorphism $\CRing
\to \Endo(\id_\catc)$.

\begin{definition}
\begin{itemize}
\item[(a)]
An $\CRing$-linear stack, that we will call $\CRing$-stack for short, is a stack $\stks$ such that $\stks(U)$ is an
$\CRing$-category for every open subset $U$, and whose restrictions
are $\CRing$-functors.  An $\CRing$-functor of $\CRing$-stacks is a
functor which is linear at the level of morphisms.  No additional
requirements are imposed on transformations of $\CRing$-functors.
\item[(b)]
Let $\cring$ be a sheaf of commutative rings on $X$.  An
$\cring$-linear stack, that we will call $\cring$-stack for short, is a $\Z$-stack $\stks$ whose center 
$\shEndo(\id_{\stks})$ is a sheaf of commutative $\cring$-algebras\footnote{By 
definition, this means that there is a morphism of sheaves of rings 
$\mu\colon \cring \to \shEndo(\id_{\stks})$. Note that the data of 
$\mu$ is equivalent to the
requirement that for every open subset $U\subset X$, and any
$\shf,\shg\in\stks(U)$ the sheaf $\shHom[\stks|_U](\shf,\shg)$ has a
structure of $\cring|_U$-module compatible with restrictions, and such
that composition is $\cring$-bilinear.}.  There is a natural notion of
$\cring$-functor\footnote{If $\stks' = (\stks',\mu')$ is another
$\cring$-stack, an $\cring$-functor $\varphi\colon \stks \to \stks'$
is a functor of $\Z$-stacks such that $\varphi(\mu(r)(\shf)) =
\mu'(r)(\varphi(\shf))$, as endomorphisms of $\varphi(\shf)$, for any $U\subset X$, $r\in\cring(U)$, and
$\shf\in\Ob(\stks(U))$.}, and transformations of $\cring$-functors are
just transformations of the underlying functors.
\end{itemize}
\end{definition}

One says that an $\cring$-functor $\varphi\colon \stks \to \stkt$ is
an equivalence (resp.\ admits a right or a left adjoint) if it is so
forgetting the $\cring$-linear structure. Note that a
quasi-inverse to $\varphi$ (resp.\ its right or left adjoint) is
necessarily an $\cring$-functor itself.  One says that
$\varphi\colon \stks \to \stkMod(\cring)$ is representable if there is
an invertible transformation $\varphi \Rightarrow
\shHom[\stks](\shf,\cdot)$ for some $\shf\in\stks(X)$.

\subsection{Operations}

Let $\stks$ and $\stks'$ be $\cring$-stacks. 
The stack $\stkHom[\cring](\stks, \stks')$ of $\cring$-functors and
transformations is an $\cring$-stack. The product 
$\stks\tens[\cring]\stks'$ is the stack associated with the
prestack $\stks \mathbin{{\mathop{\otimes}\limits^\sim}_{\raise1.5ex\hbox to-.1em{}\cring}}  \stks'$ 
defined as follows. At the level of 
objects, $\Ob(\stks \mathbin{{\mathop{\otimes}\limits^\sim}_{\raise1.5ex\hbox to-.1em{}\cring}}  \stks') = 
\Ob(\stks) \times \Ob(\stks')$. At the level of 
morphisms,
$$
\Hom[\stks \mathbin{{\mathop{\otimes}\limits^\sim}_{\raise1.5ex\hbox to-.1em{}\cring}}
\stks'(U)]((\shf_{1},\shf'_{1}),(\shf_{2},\shf'_{2})) =
\Hom[\stks(U)](\shf_{1},\shf_{2})
\tens[\cring]
\Hom[\stks'(U)](\shf'_{1},\shf'_{2}).
$$
If $\stks''$ is another $\cring$-stack, there is a natural 
$\cring$-equivalence
$$
\stkHom[\cring](\stks\tens[\cring]\stks',\stks'') \approxto
\stkHom[\cring](\stks, \stkHom[\cring](\stks',\stks'')).
$$
Let $f\colon Y\to X$ be a continuous map of topological spaces, 
$\stks$ an $\cring$-stack on $X$, and $\stkt$ an $\opb
f\cring$-stack on $Y$. Then $\oim f \stkt$ is an $\cring$-stack, 
$\opb f \stks$ is an $\opb f\cring$-stack, and there is a natural 
equivalence
\begin{equation}
\label{eq:dirinvadj}
\oim f \stkHom[\opb f\cring](\opb f\stks, \stkt) \approxto
\stkHom[\cring](\stks, \oim f \stkt).
\end{equation}

\subsection{Stacks of twisted modules}

Let $X$ be a topological space,
$\cring$ a sheaf of commutative rings on 
$X$, and $\ring$ a sheaf of not necessarily commutative
$\cring$-algebras.
Let $\catMod(\ring)$ be the category of
$\ring$-modules and $\ring$-linear morphisms. Unless otherwise stated, by $\ring$-module we mean here left 
$\ring$-module.
The prestack $\stkMod(\ring)$ of $\ring$-modules on $X$ is defined by 
$U\mapsto\catMod(\ring|_U)$, with natural restriction functors. It is 
clearly an $\cring$-stack. 

\begin{definition}
\label{de:twisted}
\begin{itemize}
\item[(a)]
A stack of $\cring$-twisted modules is an $\cring$-stack which is
locally $\cring$-equivalent to stacks of modules over
$\cring$-algebras. More precisely, an $\cring$-stack $\stkm$ is a 
stack of $\cring$-twisted modules if there exist an open covering 
$\{U_i\}_{i\in I}$ of $X$,
$\cring|_{U_i}$-algebras $\ring_i$ on $U_i$, and
$\cring|_{U_i}$-equivalences of $\cring|_{U_i}$-stacks
$\varphi_i\colon \stkm|_{U_i} \to \stkMod(\ring_i)$.
\item[(b)]
A stack of $\cring$-twisted
$\ring$-modules is an $\cring$-stack which is locally
$\cring$-equivalent to $\stkMod(\ring)$.
\item[(c)]
A stack of twisted $\cring$-modules is a stack of $\cring$-twisted
$\cring$-modules.
\end{itemize}
If $\stkm$ is a stack of $\cring$-twisted modules (resp.\ a stack of
$\cring$-twisted $\ring$-modules, resp.\ a stack of twisted
$\cring$-modules), objects of $\stkm(X)$ are called $\cring$-twisted
modules (resp.\ $\cring$-twisted $\ring$-modules, resp.\ twisted
$\cring$-modules). 
\end{definition}

Recall that a stack $\stkm$ is called additive if the
categories $\stkm(U)$ and the restriction functors are additive. 
A stack $\stkm$ is called abelian if the
categories $\stkm(U)$ are abelian, and the restriction functors are 
exact. 
Since stacks of modules over $\cring$-algebras are abelian, stacks of
$\cring$-twisted modules are also abelian. 

\begin{remark}
The stacks constructed 
in~\cite{Kashiwara1996,Kontsevich2001,Polesello-Schapira} provide 
examples of stacks of twisted modules which are of an intermediate 
nature between (a) and (b) of Definition~\ref{de:twisted}. With 
notations as in (a), denote by $\psi_i$ a quasi-inverse to $\varphi_i$. These are stacks of $\cring$-twisted modules for 
which the equivalences $\varphi_i\circ\psi_j|_{U_{ij}}$ are 
induced by isomorphisms of $\cring|_{U_{ij}}$-algebras 
$\ring_i|_{U_{ij}}\isoto\ring_j|_{U_{ij}}$. This is related to 
non-abelian cohomology as in~\cite{Giraud1971}, and we will
discuss these matters in~\cite{D'Agnolo-Polesello}.
\end{remark}

Recall that $\ring$-modules are sheaves of $\cring$-modules $\shf$ 
endowed with a morphism of sheaves of rings $m\colon\ring\to 
\shEndo[\cring](\shf)$.

\begin{definition}
If $\ring$ is an $\cring$-algebra and $\stks$ an $\cring$-stack, we denote by $\stkMod(\ring;\stks)$
the $\cring$-stack whose objects on an open subset $U\subset X$ are
pairs of an object $\shf\in\stks(U)$ and a morphism of
$\cring|_U$-algebras $m\colon \ring|_U\to\shEndo[\stks|_U](\shf)$, and
whose morphisms are those morphisms in $\stks(U)$ commuting with $m$.
We denote by $\catMod(\ring;\stks)$ the category 
$\stkMod(\ring;\stks)(X)$.
\end{definition}

Let $\ring$ and $\ringi$ be $\cring$-algebras.  Recall that an 
$\ring\tens[\cring]\ringi$-module is the same as a
$\ringi$-module $\shm$ endowed with an $\cring$-algebra morphism
$\ring\to\shEndo[\ringi](\shm)$.  Hence, there is an
$\cring$-equivalence
\begin{equation}
\label{eq:AM}
\stkMod(\ring;\stkMod(\ringi)) \approx
\stkMod(\ring\tens[\cring]\ringi).
\end{equation}
In particular, if $\stkm$ is a stack of $\cring$-twisted modules 
(resp.\ of twisted $\cring$-modules),
then $\stkMod(\ring;\stkm)$ is a stack of $\cring$-twisted
modules (resp.\ of $\cring$-twisted $\ring$-modules).

\section{Operations}\label{se:operations}

Using Morita theory, we develop the formalism of operations for 
stacks of twisted modules. We then obtain Grothendieck's six operations for derived categories of twisted modules over locally compact Hausdorff topological spaces.

\subsection{Morita theory I. Functors admitting an adjoint}\label{sse:morita1}

Morita theory describes in terms of bimodules functors between categories of modules which admit an adjoint (references are made to~\cite{Bass1968,Faith1981}). We are interested in the local analogue of this result, dealing with stacks of modules over sheaves rings. Our reference was~\cite{Kashiwara-Schapira}, where only the case of equivalences is discussed. We thus adapt here their arguments in order to deal with  functors admitting an adjoint.

\medskip 
Let $\cring$ be a sheaf of commutative rings on a topological space 
$X$, and let $\ring$ be a sheaf of not necessarily commutative
$\cring$-algebras. Denote by $\ring^\op$ the opposite algebra to 
$\ring$, given by $\ring^\op = \{a^\op \colon a\in\ring\}$ with 
product $a^\op b^\op = (ba)^\op$.  Note that left (resp.\ right) 
$\ring^\op$-modules are but right (resp.\ left) $\ring$-modules. 

For $\stks$ and $\stks'$ two $\cring$-stacks, denote by
$$
\stkHom[\cring]^\rightad(\stks,\stks')
$$ 
the full $\cring$-substack of $\stkHom[\cring](\stks,\stks')$ of
functors that admit a right adjoint. This is equivalent to the opposite of the stack of $\cring$-functors from $\stks'$ to $\stks$ that admit a left adjoint.

\begin{proposition}
\label{pr:weakMorita}
Let $\ring$ and $\ringi$ be $\cring$-algebras.  The functor
$$
\Phi\colon \stkMod(\ring\tens[\cring]\ringi^\op) \to
\stkHom[\cring]^\rightad(\stkMod(\ringi), \stkMod(\ring))
$$
given by $\shl \mapsto \shl \tens[\ringi] (\cdot)$ is an
$\cring$-equivalence.
\end{proposition}

It follows that $\cring$-functors $\stkMod(\ring) \to \stkMod(\ringi)$ which admit a left adjoint are of the form $\shHom[\ring](\shl,\cdot)$, for an $\ring\tens[\cring]\ringi^\op$-module $\shl$.

\begin{proof}
(We follow here arguments similar to those in the proof of
Morita theorem given in~\cite{Kashiwara-Schapira}.)
One checks that $\Phi$ is fully faithful.  Let us show that it is
essentially surjective.  Let $\varphi\colon \stkMod(\ringi) \to
\stkMod(\ring)$ be an $\cring$-functor admitting a right adjoint.  The
$\ring$-module $\shl = \varphi(\ringi)$ inherits a compatible
$\ringi^\op$-module structure by that of $\ringi$ itself, and we set
$\varphi'(\cdot) = \shl \tens[\ringi] (\cdot)$.  A transformation
$\alpha\colon \varphi' \Rightarrow \varphi$ is defined as follows. 
For $U\subset X$ and $\shn\in \catMod(\ringi|_U)$, the morphism
$$
\alpha(\shn) \colon \varphi(\ringi)|_U \tens[\ringi|_U] \shn \to
\varphi(\shn)
$$
is given by $l\otimes n \mapsto \varphi(\widetilde n)(l)$, where
$\widetilde n \colon \ringi|_U \to \shn$ denotes the map $b \mapsto
bn$.  We have to prove that $\alpha(\shn)$ is an isomorphism.  The 
$\ringi|_U$-module $\shn$
admits a presentation $\DSum\nolimits_j \ringi_{U_j} \to
\DSum\nolimits_i \ringi_{U_i} \to \shn \to 0$, where one sets
$(\cdot)_U = \eim u \opb u$ for $u\colon U \to X$ the open inclusion. 
We may then assume that $\shn = \DSum\nolimits_i \ringi_{U_i}$.  Since
$\varphi$ and $\varphi'$ admit a right adjoint, one has
$\varphi(\DSum\nolimits_i \ringi_{U_i}) \simeq \DSum\nolimits_i
\varphi(\ringi)_{U_i}$, and $\varphi'(\DSum\nolimits_i \ringi_{U_i}) 
\simeq
\DSum\nolimits_i \varphi'(\ringi)_{U_i}$ by 
Lemma~\ref{le:StkFunComm}.  Hence we are reduced to
prove the isomorphism $\varphi(\ringi|_U) \tens[\ringi|_U] \ringi_U
\isoto \varphi(\ringi|_U)$, which is obvious.
\end{proof}

\begin{lemma}
\label{le:StkFunComm}
Let $\ring$ and $\ringi$ be $\cring$-algebras, and let 
$\varphi\colon\stkMod(\ringi) \to \stkMod(\ring)$ be an $\cring$-functor admitting a right adjoint. Then for any family of open subsets $\{U_i\}_{i\in 
I}$ of $U\subset X$, and any $\shn\in\catMod(\ringi\vert_U)$, one has
$\varphi(\DSum\nolimits_i \shn_{U_i}) \simeq \DSum\nolimits_i
\varphi(\shn)_{U_i}$.
\end{lemma}

\begin{proof}
The proof is straightforward. We leave it to the reader to check that a functor admitting a right adjoint commutes 
with inductive limits, and in particular with direct sums. Let us check that $\varphi(\shn_{V}) \simeq \varphi(\shn)_{V}$ 
for an open inclusion $v\colon V \to U$.
Let $\psi$ be a right adjoint to $\varphi$. 
Note that the proper direct
image $\eim v$ is left adjoint to the restriction functor $\opb 
v(\cdot) = (\cdot)\vert_V$. For every $\shm\in\catMod(\ring\vert_U)$
one has
\begin{equation*}
\begin{split}
\Hom[\ring|_U](\varphi(\shn_V),\shm) &= 
\Hom[\ring|_U](\varphi(\eim v(\shn|_V)), \shm) \\
&\simeq \Hom[\ringi|_V](\shn|_V, \psi(\shm)|_V) \\
&\simeq \Hom[\ringi|_V](\shn|_V, \psi (\shm|_V) ) \\
&\simeq \Hom[\ring|_U](\eim v(\varphi(\shn|_V)), \shm) \\
&\simeq \Hom[\ring|_U](\eim v(\varphi(\shn)|_V), \shm) \\
&= \Hom[\ring|_U](\varphi(\shn)_V, \shm ),
\end{split}
\end{equation*}
where the second and fourth isomorphisms follow from the fact that 
$\psi$ and $\varphi$, respectively, are functors of stacks.
\end{proof}

\subsection{Internal product of stacks of twisted 
modules}\label{sse:inttwistedoperations}

We are now ready to define duality and internal product for stacks of twisted modules.

\medskip
Let $\cring$ be a sheaf of commutative rings on a topological space $X$.
Recall that for $\stks$ and $\stks'$ two $\cring$-stacks, we denote 
by $\stkHom[\cring]^\rightad(\stks,\stks')$ 
the stack of $\cring$-functors that admit a right adjoint.

\begin{definition}
Let $\stks$ and $\stks'$ be $\cring$-stacks on $X$.  Set
\begin{eqnarray*}
\stkinv\stks &=& \stkHom[\cring]^\rightad(\stks,\stkMod(\cring)), \\
\stks \stktens \stks' &=& 
\stkHom[\cring]^\rightad(\stkinv\stks,\stks').
\end{eqnarray*}
\end{definition}

\begin{remark}
\label{re:stktensR}
The definition of $\stkinv\stks$ does depend on the ring $\cring$, but we do not keep track of this dependence in the notation to avoid more cumbersome notations like $\stks^{\stktimes_\cring -1}$.
\end{remark}

As a consequence of Proposition~\ref{pr:weakMorita} and equivalence 
\eqref{eq:AM}, we have

\begin{proposition}
\label{pr:stkstktens}
If $\ring$ and $\ring'$ are $\cring$-algebras, there are
$\cring$-equivalences
\begin{eqnarray}
\label{eq:Astkinv}
\stkinv{\stkMod(\ring)} &\approx& \stkMod(\ring^\op), \\
\label{eq:Astktens}
\stkMod(\ring) \stktens \stkMod(\ring') &\approx&
\stkMod(\ring\tens[\cring]\ring').\\
\nonumber
&\approx&
\stkMod(\ring;\stkMod(\ring')).
\end{eqnarray}
In particular, if $\stkm$ and $\stkm'$ are stacks of $\cring$-twisted
modules on $X$, then $\stkinv\stkm$ and $\stkm \stktens \stkm'$ are
stacks of $\cring$-twisted modules on $X$.
\end{proposition}

Let us list some properties of these operations.

\begin{lemma}
    \label{le:stktens}
Let $\ring$ be an $\cring$-algebra, and $\stks$ and $\stks'$ be
$\cring$-stacks. Then there are natural $\cring$-functors
\begin{eqnarray}
\label{eq:Asigma}
\stkMod(\ring)\stktens\stks &\to& \stkMod(\ring;\stks), \\
\label{eq:sigmabidual}
\stks &\to& \stkinv{(\stkinv\stks)} = \stks \stktens \stkMod(\cring), 
\\
\label{eq:ssinvert}
\stkHom[\cring]^\rightad(\stks,\stks') &\to&
\stkHom[\cring]^\rightad(\stkinv{\stks'{}},\stkinv\stks)
= \stks' \stktens \stkinv\stks.
\end{eqnarray}
\end{lemma}

\begin{proof}
In the identification $\stkMod(\ring)\stktens\stks \approx
\stkHom[\cring]^\rightad(\stkMod(\ring^\op),\stks)$, the functor \eqref{eq:Asigma}
is given by $\phi\mapsto (\shf,m)$, where $\shf = \phi(\ring^\op)$,
and $m\colon \ring \simeq
\shEndo[\ring^\op](\ring^\op)\to[\phi]\shEndo[\stks](\phi(\ring^\op))$.
 
The functor \eqref{eq:sigmabidual} is given by $\shf \mapsto
(\phi \mapsto \phi(\shf))$, using the identification
$\stkinv{(\stkinv\stks)} =
\stkHom[\cring]^\rightad(\stkHom[\cring]^\rightad(\stks,\stkMod(\cring)),
\stkMod(\cring))$.  

Finally, the functor \eqref{eq:ssinvert} is given
by $\varphi \mapsto (\psi \mapsto \psi \circ \varphi)$, using the
identification
$\stkHom[\cring]^\rightad(\stkinv{\stks'{}},\stkinv\stks) =
\stkHom[\cring]^\rightad(\stkHom[\cring]^\rightad(\stks',\stkMod(\cring)
),\stkHom[\cring]^\rightad(\stks,\stkMod(\cring)))$.
\end{proof}

We need the following lemma from~\cite{Kashiwara-Schapira}.

\begin{lemma}
\label{le:RMtens}
For $\stkm$ a stack of $\cring$-twisted modules, there is a
natural $\cring$-functor
$$
\tens[\cring]\colon \stkMod(\cring) \times \stkm \to \stkm.
$$
\end{lemma}

\begin{proof}
For $\shm\in\catMod(\cring)$ and $\shf\in\stkm(X)$, the functor
$$
\shHom[\cring](\shm,\shHom[\stkm](\shf,\cdot)) \colon \stkm \to
\stkMod(\cring)
$$
is locally (and hence globally) representable, and we denote by
$\shm\tens[\cring]\shf$ a representative.  
\end{proof}

\begin{proposition}
\label{pr:Mhomadj}
Let $\stkm$, $\stkm'$, and $\stkm''$ be stacks of $\cring$-twisted
modules. Then there is a natural $\cring$-equivalence
$$
\stkHom[\cring]^\rightad(\stkm\stktens\stkm', \stkm'')
\approxto
\stkHom[\cring]^\rightad(\stkm, 
\stkHom[\cring]^\rightad(\stkm',\stkm'')).
$$
\end{proposition}

\begin{proof}
The above $\cring$-functor is given by $\varphi \mapsto (\shf \mapsto
(\shf' \mapsto \varphi(\psi_{\shf,\shf'})))$, where
$\psi_{\shf,\shf'}\in \stkm\stktens\stkm' =
\stkHom[\cring]^\rightad(\stkHom[\cring]^\rightad(\stkm , 
\stkMod(\cring)),
\stkm')$ is defined by $\psi_{\shf,\shf'}(\eta) = \eta(\shf)
\tens[\cring] \shf'$.  Here we used the $\cring$-functor
$\tens[\cring]$ described in Lemma~\ref{le:RMtens}.  We are then left
to prove that this functor is a local equivalence.  We may then assume
that $\stkm \approx \stkMod(\ring)$, $\stkm' \approx \stkMod(\ring')$,
and $\stkm'' \approx \stkMod(\ring'')$ for some $\cring$-algebras 
$\ring$, $\ring'$, and $\ring''$.  In this case both terms are
equivalent to $\stkMod(\ring^\op \tens[\cring] \ring^{\prime\op}
\tens[\cring] \ring'')$.
\end{proof}

\begin{proposition}
\label{pr:stktens}
Let $\ring$ be an $\cring$-algebra, and let $\stkm$, $\stkm'$, and 
$\stkm''$ be stacks of $\cring$-twisted
modules.  Then there are natural $\cring$-equivalences
\begin{eqnarray}
\label{eq:modAstks}
\stkMod(\ring)\stktens\stkm &\approx& \stkMod(\ring;\stkm), \\
\label{eq:stkminv}
\stkm &\approx& \stkinv{(\stkinv\stkm)} = \stkm \stktens 
\stkMod(\cring), \\
\label{eq:stkmcomm}
\stkm \stktens \stkm' &\approx& \stkm' \stktens \stkm, \\
\label{eq:stkinvtens}
\stkinv{(\stkm \stktens \stkm')}  &\approx& 
\stkinv\stkm \stktens \stkinv{\stkm'{}}, \\
\label{eq:stkasso}
(\stkm \stktens \stkm') \stktens \stkm''  &\approx& 
\stkm \stktens (\stkm' \stktens \stkm'').
\end{eqnarray}
\end{proposition}

\begin{proof}
Equivalences 
\eqref{eq:modAstks} and
\eqref{eq:stkminv}
follow by noticing that 
the functors
\eqref{eq:Asigma} and
\eqref{eq:sigmabidual} are local equivalences for $\stks = \stkm$.
The equivalence \eqref{eq:stkmcomm} follows by noticing that 
the functor \eqref{eq:ssinvert} is locally an equivalence for
$\stks = \stkinv\stkm$ and $\stks' = \stkm'$. 
The equivalence \eqref{eq:stkinvtens} follows from the chain of 
equivalences
\begin{equation*}
\begin{split}
\stkHom[\cring]^\rightad(\stkm\stktens\stkm', \stkMod(\cring))
&\approx
\stkHom[\cring]^\rightad(\stkm, \stkHom[\cring]^\rightad(\stkm', 
\stkMod(\cring))) \\
&\approx
\stkHom[\cring]^\rightad(\stkinv{(\stkinv\stkm)},\stkinv{\stkm'{}}).
\end{split}
\end{equation*}
The equivalence \eqref{eq:stkasso} follows from the chain of 
equivalences
\begin{equation*}
\begin{split}
\stkHom[\cring]^\rightad(\stkinv{(\stkm\stktens\stkm')}, \stkm'')
&\approx
\stkHom[\cring]^\rightad(\stkinv\stkm \stktens \stkinv{\stkm'{}}, 
\stkm'') \\
&\approx \stkHom[\cring]^\rightad(\stkinv\stkm,
\stkHom[\cring]^\rightad(\stkinv{\stkm'{}}, \stkm'')).
\end{split}
\end{equation*}
\end{proof}

Let us describe a couple of other functors.  There is a natural
$\cring$-functor
\begin{equation}
    \label{eq:Oss}
\stkMod(\cring) \to \stkm\stktens\stkinv\stkm,
\end{equation}
given by $\shf\mapsto\shf\tens[\cring](\cdot)$, in the identification
$\stkm\stktens\stkinv\stkm =
\stkHom[\cring]^\rightad(\stkinv\stkm,\stkinv\stkm)$.  Locally, 
$\stkm \approx \stkMod(\ring)$ for some $\cring$-algebra $\ring$, and 
the above functor coincides with  $\stkMod(\cring) 
\to\stkMod(\ring\tens[\cring]\ring^\op)$,
$\shf\mapsto\shf\tens[\cring]\ring$.  This has a right adjoint
$\stkMod(\ring\tens[\cring]\ring^\op) \to \stkMod(\cring)$,
$\shm\mapsto Z(\shm)$, where $Z(\shm) =
\shHom[{\ring\tens[\cring]\ring^\op}](\ring,\shm) = \{m\in\shm\colon
am=ma,\ \forall a\in\ring\}$.  Hence there is a right adjoint to
\eqref{eq:Oss}
\begin{equation}
\label{eq:ssO}
\stkm\stktens\stkinv\stkm \to \stkMod(\cring).
\end{equation}
Note also that the forgetful functor
$$
\stkMod(\ring;\stkm) \to \stkm
$$
has (locally, and hence globally) a right adjoint
$$
\ring\tens[\cring](\cdot) \colon \stkm \to \stkMod(\ring;\stkm).
$$

\subsection{Morita theory II. Relative case}

In order to describe the pull-back functor for stacks of twisted modules, we need the following relative versions of the results in Section~\ref{sse:morita1}.

\medskip
Let $f\colon Y \to X$ be a continuous map of topological spaces,
$\cring$ a sheaf of commutative rings on $X$,
$\stks$ an $\cring$-stack, and $\stkt$ an $\opb f\cring$-stack. 
Denote by
$$
\stkHom[\opb f\cring]^{\rightad\text{-}\oim f}(\opb
f\stks,\stkt)
$$ 
the full $\opb f\cring$-substack of $\stkHom[\opb f\cring](\opb
f\stks,\stkt)$ of functors $\psi$ whose image by \eqref{eq:dirinvadj}
belongs to $\stkHom[\cring]^\rightad(\stks,\oim f\stkt)$.  

\begin{proposition}
\label{pr:weakMorita2}
Let $f\colon Y \to X$ be a continuous map of topological spaces,
$\ringi$ an $\cring$-algebra on $X$, and $\ringii$ an $\opb f
\cring$-algebra on $Y$.  The functor
$$
\Psi\colon \stkMod(\ringii\tens[\opb f\cring]\opb f\ringi^\op) \to
\stkHom[\cring]^{\rightad\text{-}\oim f} (\opb f \stkMod(\ringi),
\stkMod(\ringii))
$$
given by $\shf \mapsto \shf \tens[\opb f\ringi] (\cdot)$ is an $\opb f
\cring$-equivalence.
\end{proposition}

\begin{proof}
The proof of this proposition is similar to that of 
Proposition~\ref{pr:weakMorita},
and we only show the essential surjectivity of $\Psi$.  Let
$\psi\colon \opb f\stkMod(\ringi) \to \stkMod(\ringii)$ be an $\opb
f\cring$-functor such that $\oim f \psi$ admits a right adjoint.  Set
$\shf = \psi(\opb f \ringi)$ and $\psi'(\cdot) = \shf \tens[\opb
f\ringi] (\cdot)$.  For $V\subset Y$ and $\shn\in\opb
f\stkMod(\ringi)(V)$, we have to check that the morphism
$$
\beta(\shn) \colon \psi(\opb f \ringi|_V) \tens[\opb f\ringi|_V] \shn
\to \psi(\shn),
$$
defined as the morphism $\alpha$ in Proposition~\ref{pr:weakMorita}, is an isomorphism.  By the definition
of pull-back for stacks, $\shn$ locally admits a presentation
$\DSum'_k\opb f \shm_k \to \shn$, where $\shm_k$ are objects of
$\stkMod(\ringi)$ and $\DSum'_k$ means that the sum is finite.  Thus
any $y\in Y$ has an open neighborhood $W\subset V$ such that there is
a presentation
$$
\DSum\nolimits'_k\opb f (\DSum\nolimits_j \ringi_{U_{jk}})|_W \to
\DSum\nolimits'_k\opb f (\DSum\nolimits_i \ringi_{U_{ik}})|_W \to
\shn|_W \to 0.
$$
Since $\oim f \psi$ admits a right adjoint, one has
$$
\psi(\DSum\nolimits'_k\opb f (\DSum\nolimits_i \ringi_{U_{ik}})|_W) =
\DSum\nolimits'_k\DSum\nolimits_i\psi(\opb f \ringi|_W)_{\opb
f(U_{ik})\cap W}.
$$
A similar formula holds for $\psi'$, since also $\oim f \psi'$ admits
a right adjoint.  Hence we are reduced to prove the isomorphism
$\psi(\opb f\ringi|_W) \tens[\opb f\ringi|_W] \opb f\ringi|_W \isoto
\psi(\opb f\ringi|_W)$, which is obvious.
\end{proof}

\subsection{Pull-back of stacks of twisted 
modules}\label{sse:exttwistedoperations}

We can now define the pull back of stacks of twisted modules.

\medskip
Let $f\colon Y \to X$ be a continuous map of topological spaces,
$\cring$ a sheaf of commutative rings on $X$,
$\stks$ an $\cring$-stack, and $\stkt$ an $\opb f\cring$-stack. 
Recall that we denote by
$\stkHom[\opb f\cring]^{\rightad\text{-}\oim f}(\opb
f\stks,\stkt)$ the $\opb f\cring$-stack of functors $\psi$ whose 
image by \eqref{eq:dirinvadj} admits a right adjoint.  

\begin{definition}
With the above notations, set
$$
\stkopb f\stks = \stkHom[\opb f\cring]^{\rightad\text{-}\oim f}(\opb
f(\stkinv\stks),\stkMod(\opb f\cring)).
$$
\end{definition}

\begin{remark}
Again, as in Remark~\ref{re:stktensR}, we prefer the notation $\stkopb f\stks$ to the more cumbersome $f^{\stktimes_\cring} \stks$.
\end{remark}

As a consequence of Proposition~\ref{pr:weakMorita2}, we have

\begin{proposition}
\label{pr:stkstkinv}
Let $f\colon Y \to X$ be a continuous map of topological spaces, and
$\ring$ an $\cring$-algebra on $X$.  Then, there is an $\opb
f\cring$-equivalence
\begin{equation}
\label{eq:Astkopb}
\stkopb f \stkMod(\ring) \approx \stkMod(\opb f\ring).
\end{equation}
In particular, if $\stkm$ a stack of $\cring$-twisted modules, then
$\stkopb f\stkm$ is a stack of $\opb f\cring$-twisted modules.
\end{proposition}

Let us list some properties of this operation.

\begin{proposition}
    \label{pr:stksopb}
If $\stks$ is an $\cring$-stack, there is a natural $\cring$-functor
$$
\opb f \colon \stks \to \oim f \stkopb f \stks.
$$
\end{proposition}

\begin{proof}
The usual sheaf-theoretical pull-back operation gives an
$\cring$-functor
$$
\opb f \colon \stkMod(\cring) \to \oim f \stkopb f \stkMod(\cring)
\approx \oim f \stkMod(\opb f \cring).
$$
The functor in the statement is then obtained as the composition
\begin{equation*}
\begin{split}
\stks
&\to[\eqref{eq:sigmabidual}] \stks \stktens \stkMod(\cring) \\
&\to[\id_\stks \stktens \opb f] \stks \stktens \oim f 
\stkMod(\opb f\cring) \\
&\approx \stkHom[\cring]^\rightad(\stkinv\stks, \oim f 
\stkMod(\opb f\cring)) \\
&\approx \oim f \stkHom[\opb
f\cring]^{\rightad\text{-}\oim f}(\opb f(\stkinv\stks), \stkMod(\opb f\cring)) \\
&\approx \oim f \stkopb f \stks.
\end{split}
\end{equation*}
\end{proof}

\begin{proposition}
\label{pr:Mfhomadj}
Let $\stkm$ be a stack of $\cring$-twisted
modules, and $\stkn$ a stack of $\opb f\cring$-twisted
modules. Then there is a natural $\cring$-equivalence
$$
\oim f \stkHom[\opb f \cring]^\rightad(\stkopb f \stkm, \stkn)
\approxto
\stkHom[\cring]^\rightad(\stkm, \oim f \stkn).
$$
\end{proposition}

\begin{proof}
The functor $\opb f\colon \stkm \to \oim f \stkopb f \stkm$ of 
Proposition~\ref{pr:stksopb} is locally the usual sheaf-theoretical 
pull-back, which has a right adjoint in the sheaf-theoretical 
push-forward. Moreover, it induces by \eqref{eq:dirinvadj} an 
$\opb f \cring$-functor
$$
\opb f \stkm \to \stkopb f \stkm.
$$
Hence we get a functor
$$
\stkHom[\opb f\cring]^{\rightad}(\stkopb f \stkm,\stkn) \to 
\stkHom[\opb
f\cring]^{\rightad\text{-}\oim f}(\opb f \stkm,\stkn).
$$
This is a local (and hence global) equivalence.
We thus have the chain of equivalences
\begin{equation*}
\begin{split}
\oim f \stkHom[\opb f \cring]^\rightad(\stkopb f \stkm, \stkn)
&\approx \oim f \stkHom[\opb
f\cring]^{\rightad\text{-}\oim f}(\opb f \stkm,\stkn) \\
&\approx \stkHom[\cring]^\rightad(\stkm, \oim f \stkn).
\end{split}
\end{equation*}
\end{proof}

\begin{proposition}
\label{pr:Mfformula}
Let $\stkm$ and $\stkm'$ be stacks of $\cring$-twisted
modules. Then there are natural $\opb f\cring$-equivalences
\begin{eqnarray}
    \label{eq:stkopbinv}
\stkopb f (\stkinv\stkm )  &\approx& 
\stkinv{(\stkopb f\stkm)}, \\
\label{eq:stkopbtens}
\stkopb f (\stkm \stktens \stkm') &\approx& 
\stkopb f \stkm \stktens[\opb f\cring] \stkopb f \stkm'.
\end{eqnarray}
\end{proposition}

\begin{proof}
The equivalence \eqref{eq:stkopbinv} follows from the chain of 
equivalences
\begin{equation*}
\begin{split}
\stkHom[\opb f\cring]^{\rightad\text{-}\oim f} (\opb
f(\stkinv{(\stkinv\stkm)}),\stkMod(\opb f\cring)) &\approx
\stkHom[\opb f\cring]^{\rightad\text{-}\oim f}(\opb f \stkm 
,\stkMod(\opb
f\cring)) \\
&\approx
\stkHom[\opb f\cring]^{\rightad}(\stkopb
f \stkm ,\stkMod(\opb f\cring)).
\end{split}
\end{equation*}
To prove \eqref{eq:stkopbtens}, note that, by functoriality of
$\stkopb f$, to any $\cring$-stacks $\stks$ and $\stks'$ is associated
an $\cring$-functor
$$
\stkopb f \colon \stkHom[\cring]^{\rightad} (\stks, \stks') \to \oim f
\stkHom[\opb f\cring]^{\rightad} (\stkopb f\stks, \stkopb f\stks').
$$
For $\stks = \stkinv\stkm$ and $\stks' = \stkm'$ this is locally the 
sheaf-theoretical pull-back functor 
$$
\opb f \colon \stkMod(\ring \tens[\cring] \ring')
\to \oim f \stkMod(\opb f(\ring \tens[\cring] 
\ring')), 
$$
which has a right adjoint.  Hence $\stkopb f$ has a right adjoint, 
i.e.
$$
\stkopb f \in \stkHom[\cring]^{\rightad}\bigl(
\stkHom[\cring]^{\rightad} (\stkinv\stkm, \stkm') , \oim f
\stkHom[\opb f\cring]^{\rightad} (\stkopb f(\stkinv \stkm), \stkopb f\stkm') 
\bigr).
$$
By Proposition~\ref{pr:Mfhomadj} we get a functor
$$
\stkopb f \stkHom[\cring]^{\rightad} (\stkinv\stkm, \stkm') \to
\stkHom[\opb f\cring]^{\rightad} (\stkopb f(\stkinv\stkm), \stkopb
f\stkm').
$$
This is locally, and hence globally, an equivalence.
\end{proof}

\subsection{Twisted sheaf-theoretical operations} 

Let us now show how the usual operations of sheaf theory extend to
the twisted case. For the classical non-twisted case, that we do not 
recall here,  we refer e.g.\ to
\cite{Kashiwara-Schapira1990}.

\begin{proposition}
\label{pr:operations}
Let $f\colon Y\to X$ be a continuous map of topological spaces, and
$\stkm$ and $\stkm'$ be stacks of $\cring$-twisted modules.  Then
there exist $\cring$-functors
\begin{align*}
\tens[\cring] & \colon \stkm \times \stkm' \to \stkm \stktens \stkm',
\\
\shHom[\cring] & \colon (\stkinv\stkm)^\op \times \stkm' \to \stkm
\stktens \stkm', \\
\opb f & \colon \stkm \to \oim f \stkopb f \stkm, \\
\oim f & \colon \oim f \stkopb f \stkm \to \stkm.
\end{align*}
If moreover $X$ and $Y$ are locally compact Hausdorff topological spaces, there exists an $\cring$-functor
$$
\eim f \colon \oim f \stkopb f \stkm \to \stkm.
$$
If $U\subset X$ is an open subset where $\stkm|_U \approx
\stkMod(\ring)$ and $\stkm'|_U \approx \stkMod(\ring')$ for some
$\cring|_U$-algebras $\ring$ and $\ring'$, then the restrictions to
$U$ of the above functors coincide with the usual sheaf operations.
\end{proposition}

\begin{proof}
The functor
$$
\tens[\cring]\colon \stkm \times \stkm' \to \stkm \stktens
\stkm' =
\stkHom[\cring]^\rightad(\stkHom[\cring]^\rightad(\stkm,\stkMod(\cring)),\stkm')
$$
is defined by $(\shf,\shf')\mapsto(\phi\mapsto
\phi(\shf)\tens[\cring]\shf')$, using Lemma~\ref{le:RMtens}.

For $\shf$ an object of $\stkinv\stkm$ there is a natural functor
\begin{equation}
\label{eq:homOtemp}
\stkm \stktens \stkm' \to \stkm'
\end{equation}
given by $\phi \mapsto \phi(\shf)$ in the identification $\stkm
\stktens \stkm' = \stkHom[\cring]^\rightad(\stkinv\stkm,\stkm')$.  
Locally
this corresponds to the functor $\stkMod(\ring\tens[\cring]\ring') \to
\stkMod(\ring')$, $\shm\mapsto \shf\tens[\ring]\shm$, for 
$\shf\in\catMod(\ring^\op)$.  If $\shn$ is an
$\ring'$-module, there is a functorial isomorphism
$\shHom[\ring'](\shf\tens[\ring]\shm,\shn) \simeq
\shHom[{\ring\tens[\cring]\ring'}](\shm,\shHom[\cring](\shf,\shn))$. 
Hence \eqref{eq:homOtemp} admits a right adjoint, that we denote by
$\shHom[\cring](\shf,\cdot)$.  This construction is functorial in
$\shf$, and hence we get the bifunctor $\shHom[\cring](\cdot,\cdot)$.

The functor $\opb f$ was constructed in Proposition~\ref{pr:stksopb}.

The functor $\oim f$ is obtained by noticing that if $\stkm$ is a
stack of $\cring$-twisted modules, then $\opb f$ is locally the usual
sheaf-theoretical pull-back, which admits a right adjoint.

Assume that $f\colon Y \to X$ is a continuous map of locally
compact Hausdorff topological spaces.  Recall that for an $\opb 
f\ring$-module $\shg$ on $Y$
one denotes by $\eim f\shg$ the subsheaf of $\oim f \shg$ of sections
$s\in\oim f \shg(U)$ such that $f|_{\supp(s)}$ is proper.  Such a
condition is local on $X$, and hence for a stack of
$\cring$-twisted modules $\stkm$ there is an $\cring$-functor
$\eim f \colon \oim f \stkopb f \stkm \to \stkm$
locally given by the usual proper direct image functor for sheaves
just recalled.
\end{proof}

\subsection{Derived twisted operations} 
\label{sse:twistedgrothendieck}

Let us now deal with the twisted version of Grothendieck's formalism of six operations for 
sheaves over locally compact Hausdorff topological spaces.  We do not 
recall here such formalism for the
classical non-twisted case, referring instead e.g.\ to
\cite{Kashiwara-Schapira1990}.

\medskip

Let $\stkm$ be a stack of $\cring$-twisted modules, and denote by 
$\TDC(\stkm)$ the derived category of the abelian category 
$\stkm(X)$.  Let
$\BDC(\stkm)$ (resp.\ $\PDC(\stkm)$, resp.\ $\NDC(\stkm)$) be the full
triangulated subcategory of $\TDC(\stkm)$ whose objects have bounded
(resp.\ bounded below, resp.\ bounded above) amplitude.  

\begin{lemma}
The category $\stkm(X)$ has enough injective objects.
\end{lemma}

\begin{proof}
The classical proof, found e.g.\ in~\cite[Proposition~2.4.3]{Kashiwara-Schapira1990}, 
adapts as follows. Consider the natural map $p\colon \hat X \to X$, 
where $\hat X$ is the set $X$ endowed with the discrete topology. 
For $F\in \stkm(X)$, the adjunction morphism $F\to \oim p \opb p 
F$ is injective, and the functor $\oim p$ is left exact. It is thus 
enough to find an injection $\opb p F \to I$, where $I$ is an 
injective object in $\stkopb p \stkm(\hat X)$. Since $\hat X$ is 
discrete, $\stkopb p \stkm$ is equivalent to a stack of (non 
twisted) modules.\footnote{Another proof is obtained by applying Grothendieck's criterion, stating that a category has enough injective objects if it admits small filtrant inductive limits, which are exact, and if it admits a generator. Let 
$\{U_i\}_{i\in I}$ be an open covering of $X$, let $\varphi_i\colon 
\stkMod(\ring_i) \approxto \stkm|_{U_i}$ be $\cring$-equivalences for 
some $\cring$-algebras $\ring_i$, and let $G_i$ be generators of 
$\catMod(\ring_i)$. Then a generator of $\stkm(X)$ is given by 
$G=\DSum_i j_{i!} \varphi_i(G_i)$, where $j_i\colon U_i \to X$ are 
the open inclusions.}
\end{proof}

Let $f\colon Y\to X$ be a continuous map of topological spaces. 
Deriving the functors $\opb f$, $\oim f$, and $\hom[\cring]$, one 
gets functors
\begin{align*}
\opb f & \colon \TDC^{\pm,\mathrm{b}}(\stkm) \to 
\TDC^{\pm,\mathrm{b}} (\stkopb f \stkm), \\
\roim f & \colon \PDC( \stkopb f \stkm) \to \PDC(\stkm), \\
\rhom[\cring] & \colon \NDC(\stkinv\stkm)^\op \times \PDC(\stkm') \to
\PDC(\stkm \stktens \stkm').
\end{align*}
Assuming that the weak global dimension of $\cring$ is finite, one 
gets that
$\stkm(X)$ has enough flat objects. Deriving $\tens[\cring]$ one gets 
a functor
$$
\ltens[\cring] \colon \TDC^{\pm,\mathrm{b}}(\stkm) \times
\TDC^{\pm,\mathrm{b}}(\stkm') \to \TDC^{\pm,\mathrm{b}}(\stkm \stktens
\stkm').
$$
Assuming that $f$ is a map between locally compact Hausdorff 
topological spaces, one can derive the functor $\eim f$, and get
$$
\reim f \colon \PDC( \stkopb f \stkm) \to \PDC(\stkm).
$$
Assume that $\eim f$ has finite cohomological dimension.
The usual construction of Poincar\'e-Verdier duality (cf 
e.g.~\cite[\S 3.1]{Kashiwara-Schapira1990}) extends to the twisted 
case as follows.
Let $L\in \catMod(\Z_X)$, and consider the functor
$$
\eim f (\cdot \tens[\Z] L) \colon \oim f \stkopb f \stkm \to \stkm.
$$
Denote by $\shi(\stkm)$ the full substack of $\stkm$ of injective 
objects. 
Assuming that $L$ is flat and $f$-soft, there exists a functor
$$
\epb f_L \colon \shi(\stkm) \to \shi(\oim f\stkopb f \stkm)
$$
characterized by the isomorphism, functorial in $I$ and $G$,
$$
\hom[\stkm](\eim f (G \tens[\Z] L), I)
\simeq
\oim f \hom[\stkopb f \stkm](G , \epb f_L I).
$$
In fact, the above isomorphism shows that the existence of $\epb f_L$ 
is a local problem, and locally this is the classical construction. 
As in the classical case, one finally gets a functor
$$
\epb f \colon \PDC(\stkm) \to \PDC( \stkopb f \stkm)
$$
by letting $\epb f F$ be the simple complex associated to the double 
complex $\epb f_{L^\bullet}I^\bullet$, where $I^\bullet\in 
K^+(\shi(\stkm)(X))$ is quasi-isomorphic to $F$, and $L^\bullet$ is a 
(non twisted) bounded, flat, $f$-soft resolution of $\Z_Y$.

One proves that the usual formulas relating the six operations above, 
like adjunction, base-change, or projection formulas, hold.

\section{Descent}\label{se:descent}

Effective descent data for stacks of twisted modules, called twisting
data, are considered
in~\cite{Kashiwara1989,Kashiwara-Schmid1994,Kashiwara-Schapira}, and we recall here this notion using the language of semisimplicial complexes.  
We then describe in terms of twisting data equivalences, operations, and the example of twisted modules associated with a line bundle.

\subsection{Morita theory III. Equivalences}

In Section~\ref{sse:morita1} we recalled how functors between stacks of modules admitting an adjoint are described in term of bimodules. We discuss here the particular case of equivalences. (References are again made to~\cite{Bass1968,Faith1981,Kashiwara-Schapira}.)

\medskip
Two $\cring$-algebras $\ring$ and $\ringi$ are called
Morita equivalent if $\stkMod(\ring)$ and $\stkMod(\ringi)$ are
$\cring$-equivalent. Let us recall how such equivalences are 
described in terms of $\ring\tens[\cring]\ringi^{\op}$-modules.

\begin{proposition}
\label{pr:preMorita}
Let $\shl$ be an $\ring\tens[\cring]\ringi^{\op}$-module.  Then the
following conditions are equivalent:
\begin{itemize}
\item[(i)] There exists a $\ringi\tens[\cring]\ring^{\op}$-module
$\shl'$, such that $\shl \tens[\ringi] \shl' \simeq \ring$ as
$\ring\tens[\cring]\ring^{\op}$-modules and $\shl' \tens[\ring] \shl
\simeq \ringi$ as $\ringi\tens[\cring]\ringi^{\op}$-modules. 
\item[(ii)] For $\shl^{*\ring} = \shHom[\ring](\shl,\ring)$, the
canonical morphism $\shl \tens[\ringi] \shl^{*\ring}\to \ring$ is an
isomorphism of $\ring\tens[\cring]\ring^{\op}$-modules, and
$\shl^{*\ring} \tens[\ring] \shl \simeq \ringi$ as
$\ringi\tens[\cring]\ringi^{\op}$-modules.  

\item[(iii)] $\shl$ is a
faithfully flat $\ring$-module locally of finite presentation, and
there is an $\cring$-algebra isomorphism $\ringi^\op \isoto
\shEndo[\ring](\shl)$.  

\item[(iv)] $\shl$ is a faithfully flat
$\ringi^\op$-module locally of finite presentation, and there is an
$\cring$-algebra isomorphism $\ring \isoto 
\shEndo[\ringi^\op](\shl)$.

\item[(v)] $\shl\tens[\ringi] (\cdot) \colon \stkMod(\ringi) \to
\stkMod(\ring)$ is an $\cring$-equivalence.  

\item[(vi)]
$\shHom[\ring](\shl,\cdot)\colon \stkMod(\ring) \to \stkMod(\ringi)$
is an $\cring$-equivalence.
\end{itemize}
\end{proposition}

\begin{definition}
\label{de:morita}
An $\ring\tens[\cring]\ringi^{\op}$-module $\shl$ is called invertible
if the equivalent conditions in Proposition~\ref{pr:preMorita} are
satisfied.
\end{definition}

The ring $\ring$ itself is the invertible $\ring\tens[\cring]\ring^{\op}$-module corresponding to the identity functor of $\stkMod(\ring)$. Note that invertible $\ring\tens[\cring]\ring^{\op}$-modules are not necessarily locally isomorphic to $\ring$ as $\ring$-modules, even if $\ring$ is a commutative ring.

\begin{theorem}[Morita]
\label{th:Morita}
If $\varphi\colon\stkMod(\ringi) \to \stkMod(\ring)$ is an equivalence
of $\cring$-stacks, then $\shl = \varphi(\shb)$ is an invertible
$\ring\tens[\cring]\ringi^{\op}$-module, and $\varphi \simeq \shl
\tens[\ringi] (\cdot)$.  Moreover, a quasi-inverse to $\varphi$ is
given by $\hom[\ring](\shl, \cdot) \simeq \shl^{*\ring} \tens[\ring]
(\cdot)$.
\end{theorem}

\begin{proof}
Let $\psi$ be a quasi-inverse to $\varphi$.  Since $\psi$ is right
adjoint to $\varphi$, by Proposition~\ref{pr:weakMorita} $\shl =
\varphi(\shb)$ is an $\ring\tens[\cring]\ringi^{\op}$-module such that
$\varphi \simeq \shl \tens[\ringi] (\cdot)$.  Interchanging the role
of $\varphi$ and $\psi$ there also exists a
$\ringi\tens[\cring]\ring^{\op}$-module $\shl'$ such that $\psi \simeq
\shl' \tens[\ring] (\cdot)$.  Since $\varphi \circ \psi$ and $\psi
\circ \varphi$ are isomorphic to the identity functors, $\shl$ is
invertible and $\shl'\simeq \shl^{*\ring}$.  Finally, since $\shl$ is 
a flat $\ring$-module locally of finite presentation, one has 
$\shl^{*\ring} \tens[\ring] (\cdot) \simeq \shHom[\ring](\shl,\cdot)$.
\end{proof}

\subsection{Twisting data on an open covering}

By definition, if $\stkm$ is a stack of $\cring$-twisted modules
there exist an open covering $\{U_i\}_{i\in I}$ of $X$,
$\cring|_{U_i}$-algebras $\ring_i$ on $U_i$, and
$\cring|_{U_i}$-equivalences of $\cring|_{U_i}$-stacks
$\varphi_i\colon \stkm|_{U_i} \to \stkMod(\ring_i)$.  Let $\psi_{i}$
be a quasi-inverse of $\varphi_{i}$, and let
$\alpha_{i}\colon \psi_{i}\circ \varphi_{i} \Rightarrow \id_{\stkm}$
be an invertible transformation.  By (the $\cring$-linear analogue of)
Proposition~\ref{pr:patch}, the following descent datum for stacks is
enough to reconstruct $\stkm$
\begin{equation}
    \label{eq:isomodA} (\{U_{i}\}_{i\in I},
    \{\stkMod(\ring_{i})\}_{i\in I}, \{\varphi_{ij}\}_{i,j\in I},
    \{\alpha_{ijk}\}_{i,j,k\in I}).
\end{equation}
Here $\varphi_{ij} =
\varphi_{i}\vert_{U_{ij}}\circ\psi_{j}\vert_{U_{ij}}$, and
$\alpha_{ijk}\colon \varphi_{ij}\circ \varphi_{jk} \Rightarrow
\varphi_{ik}$ is induced by $\alpha_{j}$, so that they satisfy
condition \eqref{eq:phiijcommut}.  Functors as $\varphi_{ij}$ are
described by Morita's Theorem~\ref{th:Morita}, so that the descent 
datum
\eqref{eq:isomodA} is replaced by
\begin{equation}
     \label{eq:twistdat0} \twst = (\{U_{i}\}_{i\in I},
     \{\ring_{i}\}_{i\in I}, \{\shl_{ij}\}_{i,j\in I},
     \{a_{ijk}\}_{i,j,k\in I}),
\end{equation}
where $\{U_{i}\}_{i\in I}$ is an open covering of $X$, $\ring_i$ is an
$\cring|_{U_i}$-algebra on $U_i$, $\shl_{ij}$ are invertible
$\ring_i\tens[\cring]\ring_j^{\op}\vert_{U_{ij}}$-modules, and
$a_{ijk}\colon \shl_{ij}\tens[\ring_j] \shl_{jk}\vert_{U_{ijk}} \to
\shl_{ik}\vert_{U_{ijk}}$ are isomorphisms of
$\ring_i\tens[\cring]\ring_k^{\op}\vert_{U_{ijk}}$-modules satisfying
the analogue of condition \eqref{eq:phiijcommut}.  As in the proof of
Proposition~\ref{pr:patch}, up to equivalence a twisted module
$\shf\in\stkm(X)$ is thus described by a pair
$$
(\{\shf_{i}\}_{i\in I}, \{m_{ij}\}_{i,j\in I}),
$$
where $\shf_{i}\in\catMod(\ring_{i})$, and $m_{ij}\colon \shl_{ij}
\tens[\ring_{j}] \shf_j|_{U_{ij}} \to \shf_i|_{U_{ij}}$ is an
isomorphism of $\ring_i|_{U_{ij}}$-modules on $U_{ij}$ such that the
following diagram on $U_{ijk}$ commutes
$$
\xymatrix@C=3.5em{ \shl_{ij}\tens[\ring_j] \shl_{jk}\tens[\ring_k] 
\shf_{k}
\ar[r]^-{a_{ijk} \tens \id_{\shf_{k}}} \ar[d]^{\id_{\shl_{ij}}\tens 
m_{jk}} & \shl_{ik}\tens[\ring_k] \shf_{k}
\ar[d]^{m_{ik}} \\ \shl_{ij}\tens[\ring_j] \shf_{j} \ar[r]^-{m_{ij}} &
\shf_{i} .  }
$$
This is actually the definition of twisted modules given
in~\cite{Kashiwara1989}.  It is also an example of twisting data, of
which we now give a more general definition.

\subsection{Twisting data}\label{sse:descent}

We shall use here the language of semisimplicial complexes. On the one hand, this allows one to consider more general situations than open coverings, on the other hand, it provides a very efficient bookkeeping of indices. 

\medskip 
Recall that semisimplicial complexes are diagrams of continuous maps of
topological spaces\footnote{In dealing with stacks, we will only need the terms $X^{[r]}$ with 
$r\leq 3$.}
\begin{equation}
\label{eq:simplicial}
\xymatrix@C=10ex{ X^{[3]} \ar@<-1.5ex>[r]\ar@<-.5ex>[r]
\ar@<.5ex>[r]\ar@<1.5ex>[r]^-{q^{[3]}_0,\dots, q^{[3]}_3} & X^{[2]}
\ar@<-1ex>[r]\ar[r]\ar@<1ex>[r]^-{q^{[2]}_0, q^{[2]}_1, q^{[2]}_2} &
X^{[1]} \ar@<-.5ex>[r]\ar@<.5ex>[r]^-{q^{[1]}_0, q^{[1]}_1} & X^{[0]}
\ar[r]^-{q^{[0]}_0 = q^{[0]} = q} & X^{[-1]} = X, }
\end{equation}
satisfying the commutativity relations 
$$
q^{[r]}_j \circ q^{[r+1]}_i =
q^{[r]}_i \circ q^{[r+1]}_{j+1}, 
$$
for $0\leq i \leq j \leq r$.
In the coskeleton construction, one considers the topological space
$$
X^{[r]+1} = \{ (x_0,\dots,x_{r+1})\in (X^{[r]})^{r+2} \colon
q^{[r]}_j(x_i) = q^{[r]}_i(x_{j+1}) \text{ for } 0\leq i \leq j \leq r
\},
$$
and let $q^{[r+1]}\colon X^{[r+1]} \to X^{[r]+1}$ be the map $x
\mapsto (q^{[r+1]}_0(x),\dots,q^{[r+1]}_{r+1}(x))$.  
Hence there are commutative diagrams for $0\leq i \leq r+1$
$$
\xymatrix@C=10ex{ 
X^{[3]} \ar[r]^-{q^{[3]}_i} 
\ar[d]^{q^{[3]}} &
X^{[2]} \ar[r]^-{q^{[2]}_i} 
\ar[d]^{q^{[2]}} &
X^{[1]} \ar[r]^-{q^{[1]}_i} 
\ar[d]^{q^{[1]}} &
X^{[0]} \ar[r]^-{q} 
\ar[d]^{q^{[0]}} & X \\
X^{[2]+1} \ar[ur] &
X^{[1]+1} \ar[ur] &
X^{[0]+1} \ar[ur] &
X \ar@{=}[ur], }
$$
where the diagonal arrows are the projection to the 
$i$th factor.  

\begin{example}
\label{ex:twist}
\begin{itemize}
\item[(a)]
Let us say that a semisimplicial complex is coskeletal if $X^{[r+1]} \simeq X^{[r]+1}$ for $r\geq 0$.
In other words, $X^{[r]} = Y\times_X \cdots \times_X Y$ is the $(r+1)$-fold fibered product of a continuous map $q\colon Y \to X$,
and $q^{[r]}_i$ the projection omitting the $i$th factor.
\item[(a1)]
A particular case of coskeletal semisimplicial complex is the one attached to an open covering $\{U_{i}\}_{i\in I}$ of $X$. In this case, $Y = \DUnion\nolimits_{i\in I}U_{i}$ is the
disjoint union of the $U_{i}$'s, and $q$ is the natural map (which
is a local homeomorphism).  Note that $X^{[r]} =
\DUnion\nolimits_{i_0,\dots,i_r\in I}U_{i_0\cdots i_r}$.
\item[(a2)]
Another particular case of coskeletal semisimplicial complex is obtained when $q\colon Y\to X$ is a
principal $G$-bundle, for $G$ a topological group.  Denoting by
$m\colon G\times Y \to Y$ the group action, this semisimplicial complex
is identified with
$$
\xymatrix{ G \times G \times G \times Y \ar@<-1.5ex>[r]\ar@<-.5ex>[r]
\ar@<.5ex>[r]\ar@<1.5ex>[r]^-{q^{[3]}_i} & G \times G \times Y
\ar@<-1ex>[r]\ar[r]\ar@<1ex>[r]^-{q^{[2]}_i} & G \times Y
\ar@<-.5ex>[r]\ar@<.5ex>[r]^-{q^{[1]}_i} & Y \ar[r]^-{q} & X, }
$$
where $q^{[r]}_r = \id_{G^{r-1}} \times m$, $q^{[r]}_0$ is the 
projection omitting the $0$th factor, and 
$q^{[r]}_i(g_0,\dots,g_{r-1},y) = 
(g_0,\dots,g_{i-1}g_i,\dots,g_{r-1},y)$ for $0\leq i < r$.
\item[(b)]
Other examples of semisimplicial complexes are the ones attached to
hypercoverings, where $X^{[r+1]}$ is induced by an open covering of $X^{[r]+1}$. These are of the form $X^{[0]} =
\DUnion\nolimits_{i}U_{i}$ for
$X=\Union\nolimits_{i\in I} U_i$, 
$X^{[1]} = \DUnion\nolimits_{i,j,\alpha}U_{ij}^\alpha$,
for
$U_{ij} = \Union\nolimits_{\alpha\in A_{ij}} U_{ij}^\alpha$,
$X^{[2]} = \DUnion\nolimits_{i,j,k,\alpha,\beta,\gamma,\xi} U^\xi_{ijk\alpha\beta\gamma}$
for
$U_{ij}^\alpha \cap
U_{kj}^\beta \cap U_{ki}^\gamma = \Union\nolimits_{\xi\in \Xi_{ijk}^{\alpha\beta\gamma}}U^\xi_{ijk\alpha\beta\gamma}$, 
and so on.
\end{itemize}
\end{example}

Let $s > r$, $0\leq i_0 < \cdots < i_r \leq s$, and $0\leq i_{r+1} <
\cdots < i_s \leq s$, be such that $\{i_0,\dots, i_s\} =
\{0,\dots,s\}$.  If $\shf$ is a sheaf on $X^{[r]}$, we denote by
$\shf_{i_0 \cdots i_r} = \opb{(q^{[r+1]}_{i_{r+1}}\circ \cdots \circ
q^{[s]}_{i_s})} \shf$ its sheaf-theoretical pull-back to
$X^{[s]}$, and we use the same notations for
morphisms of sheaves.\footnote{In the coskeletal case,
$\shf_{i_0\cdots i_{r}}$ is the pull-back of $\shf$ by the projection
to the $(i_0,\dots,i_r)$th factors.}

\begin{definition}
\begin{itemize}
\item[(i)]
An $\cring$-twisting datum on $X$ is a quadruplet\footnote{This notion
was discussed in \cite{Kashiwara1989} for semisimplicial complexes
attached to open coverings, and in
\cite{Kashiwara-Schmid1994} for coskeletal semisimplicial complexes.}
\begin{equation}
     \label{eq:twistcov} \twst = (X^{[\bullet]} \to[q] X, \ring, \shl,
     a),
\end{equation}
where $X^{[\bullet]} \to[q] X$ is a semisimplicial complex, $\ring$ 
is a $\opb{q} \cring$-algebra on
$X^{[0]}$, $\shl$ is an invertible $\ring_0
\tens[\cring]\ring_1^{\op}$-module on $X^{[1]}$, and $a\colon
\shl_{01} \tens[\ring_1] \shl_{12} \to \shl_{02}$ is an isomorphism of
$\ring_0 \tens[\cring]\ring_2^{\op}$-modules on $X^{[2]}$ such that
the following diagram on $X^{[3]}$ commutes\footnote{Let us denote by
$\shl^{[r]}_{ij}$ the sheaf $\shl_{ij}$ on $X^{[r]}$.  Then one should
pay attention to the fact that in $X^{[3]}$ one has isomorphisms like
$\shl^{[3]}_{01} \simeq \opb{(q^{[3]}_{3})}\shl^{[2]}_{01}$, but not
equalities.  Thus, much as in Definition~\ref{de:stack}~(iv), one
should write \eqref{eq:La} more precisely as
$$
\xymatrix@C=2ex@R=3ex{ \opb{(q^{[3]}_3)} (\shl^{[2]}_{01}
\tens[\ring_1] \shl^{[2]}_{12} ) \tens[\ring_2] \shl^{[3]}_{23}
\ar[d]_-{a\tens \id_{\shl^{[3]}_{23}}} & \shl^{[3]}_{01} 
\tens[\ring_1] \shl^{[3]}_{12}
\tens[\ring_2] \shl^{[3]}_{23} \ar@{-}[r]^-\sim \ar@{-}[l]_-\sim &
\shl^{[3]}_{01} \tens[\ring_1] \opb{(q^{[3]}_0)} ( \shl^{[2]}_{01}
\tens[\ring_1] \shl^{[2]}_{12} ) \ar[d]_-{\id_{\shl^{[3]}_{01}} \tens 
a} \\
\opb{(q^{[3]}_3)} \shl^{[2]}_{02} \tens[\ring_2] \shl^{[3]}_{23} &&
\shl^{[3]}_{01} \tens[\ring_1] \opb{(q^{[3]}_0)} \shl^{[2]}_{02} \\
\shl^{[3]}_{02} \tens[\ring_2] \shl^{[3]}_{23} \ar@{-}[u]^-\sim
\ar@{-}[d]_-\sim && \shl^{[3]}_{01} \tens[\ring_1] \shl^{[3]}_{13}
\ar@{-}[u]^-\sim \ar@{-}[d]_-\sim \\
\opb{(q^{[3]}_1)} ( \shl^{[2]}_{01} \tens[\ring_1] \shl^{[2]}_{12} )
\ar[d]_-a && \opb{(q^{[3]}_2)} ( \shl^{[2]}_{01} \tens[\ring_1]
\shl^{[2]}_{12} ) \ar[d]_-a \\
\opb{(q^{[3]}_1)} \shl^{[2]}_{02} & \shl^{[3]}_{03} \ar@{-}[r]^-\sim
\ar@{-}[l]_-\sim & \opb{(q^{[3]}_2)} \shl^{[2]}_{02} . }
$$
Such a level of precision is both quite cumbersome and easy to attain,
so we prefer a sloppier but clearer presentation.}
\begin{equation}
\label{eq:La}
\xymatrix@C=5em{ \shl_{01}\tens[\ring_1] \shl_{12}\tens[\ring_2] \shl_{23}
\ar[r]^-{a_{012} \tens \id_{\shl_{23}}} \ar[d]^{\id_{\shl_{01}}\tens 
a_{123}} & \shl_{02}\tens[\ring_2] \shl_{23}
\ar[d]^{a_{023}} \\ \shl_{01}\tens[\ring_1] \shl_{13}
\ar[r]^-{a_{013}} & \shl_{03} .  }
\end{equation}
\item[(ii)]
A coskeletal $\cring$-twisting datum on $X$ is an $\cring$-twisting datum whose associated semisimplicial complex is coskeletal.
\end{itemize}
\end{definition}

One can now mimic the construction in the sketch of proof of
Proposition~\ref{pr:patch}.  Denote by $\catMod(\twst)$ the category
whose objects are pairs $(\shf,m)$, where $\shf$ is an $\ring$-module
on $X^{[0]}$, and $m\colon \shl\tens[\ring_{1}] \shf_1 \to \shf_0$ is
an isomorphism of $\ring_0$-modules on $X^{[1]}$ such that the
following diagram on $X^{[2]}$ commutes
\begin{equation}
\label{eq:LFm}
\xymatrix{ \shl_{01}\tens[\ring_1] \shl_{12}\tens[\ring_2] \shf_{2}
\ar[r]^-{a \tens \id_{\shf_{2}}} \ar[d]^{\id_{\shl_{01}}\tens m_{12}} 
& \shl_{02}\tens[\ring_2] \shf_{2}
\ar[d]^{m_{02}} \\ \shl_{01}\tens[\ring_1] \shf_{1} \ar[r]^-{m_{01}} &
\shf_{0} , }
\end{equation}
and whose morphisms $\alpha \colon (\shf,m) \to (\shf',m')$ consists
of morphisms of $\ring$-modules $\alpha \colon \shf \to \shf'$,
such that the following diagram on $X^{[1]}$ commutes
$$
\xymatrix{ \shl\tens[\ring_{1}] \shf_1 \ar[r]^-{m} 
\ar[d]^{\id_\shl\tens\alpha_{1}} &
\shf_0 \ar[d]^{\alpha_{0}} \\
\shl\tens[\ring_{1}] \shf'_1 \ar[r]^-{m'} & \shf'_0 .  }
$$

\begin{definition}
Let $\twst$ be an $\cring$-twisting datum on $X$.  We denote by
$\stkMod(\twst)$ the prestack on $X$ defined by $U \mapsto
\catMod(\twst\vert_{U})$, which is in fact an $\cring$-stack. Here, $\twst\vert_{U}$ denotes the $\cring\vert_{U}$-twisting datum on $U$ naturally induced by $\twst$.
\end{definition}

Note that if $\ringi$ is an $\cring$-algebra on $X$, then
$\stkMod(\ringi) \approx \stkMod(\twid_\ringi)$ for
$$
\twid_\ringi = (X\to[\id]X, \ringi, \ringi, \cdot)
$$
the trivial $\cring$-twisting datum, with $\cdot$ being the canonical
isomorphism $\ringi\tens[\ringi]\ringi\isoto\ringi$.

\medskip
We spend the rest of this section to show that $\stkMod(\twst)$ is
actually a stack of $\cring$-twisted modules, using arguments adapted from those in~\cite{Kashiwara-Schmid1994}.
In order to get this result it seems natural to assume that the maps $q^{[r]} \colon X^{[r]} \to X^{[r-1]+1}$ admit local sections for $r=0,1,2,3$. However, we will consider here the
stronger assumption
\begin{equation}
\label{eq:locsectstrict}
\text{the maps $q^{[r]}$, for $r=0,1,2,3$, admit sections locally on 
$X$.}
\end{equation}
Note that for coskeletal semisimplicial complexes this reduces to the assumption
\begin{equation}
\label{eq:cosklocsectstrict}
\text{$q \colon X^{[0]} \to X$ admits local sections,}
\end{equation}
which holds true for semisimplicial complexes attached to open coverings or to principal $G$-bundles, as in
Example~\ref{ex:twist} (a1) and (a2). In general, \eqref{eq:locsectstrict} does not hold for semisimplicial complexes attached to hypercoverings, as in Example~\ref{ex:twist} (b).

Let $\twst' = (X^{\prime [\bullet]} \to[q'] X, \ring', \shl', a')$ be
another $\cring$-twisting datum on $X$.

\begin{definition}
\begin{itemize}
\item[(i)]
A refinement of $\cring$-twisting data $\rho\colon \twst' \to \twst$
consists of commutative diagrams
$$
\xymatrix{ X^{\prime[3]} \ar[r]^{q_k^{\prime[3]}} \ar[d]^{\rho^{[3]}} 
&
X^{\prime[2]} \ar[r]^{q_j^{\prime[2]}} \ar[d]^{\rho^{[2]}} & 
X^{\prime[1]}
\ar[r]^{q_i^{\prime[1]}} \ar[d]^{\rho^{[1]}} & X^{\prime[0]} 
\ar[r]^{q'}
\ar[d]^{\rho^{[0]}} & X \ar@{=}[d] \\
X^{[3]} \ar[r]^{q_k^{[3]}} & X^{[2]} \ar[r]^{q_j^{[2]}} &
X^{[1]} \ar[r]^{q_i^{[1]}} & X^{[0]} \ar[r]^{q} & X, }
$$
of an isomorphism of $\opb q \cring$-algebras $\opb{(\rho^{[0]})}\ring
\isoto \ring'$, and of an isomorphism of
$\ring'_0\tens[\cring]\ring^{\prime\op}_1$-modules
$\opb{(\rho^{[1]})}\shl \isoto \shl'$ compatible with $a$ and $a'$.
\item[(ii)]
To $\rho\colon \twst' \to \twst$ one associates the functor $\opb \rho
\colon \stkMod(\twst) \to \stkMod(\twst')$, given by $(\shf,m) \mapsto
(\opb{(\rho^{[0]})}\shf,\opb{(\rho^{[1]})}m)$.
\end{itemize}
\end{definition}

\begin{lemma}
\label{le:globsect}
Let $\twst = (X^{[\bullet]} \to[q] X, \ring, \shl, a)$ be such that
the maps $q^{[r]}$ admit global sections $s^{[r]}$. Then
\begin{itemize}
\item[(i)]
there is a
refinement of $\cring$-twisting data $\tilde s\colon
\twid_{\opb{(s^{[0]})}\ring} \to \twst$,
\item[(ii)]
the functor $\opb{\tilde s}\colon \stkMod(\twst) \to \stkMod(\opb{(s^{[0]})}\ring)$ is an equivalence.
\end{itemize}
\end{lemma}

\begin{proof}
Define the maps $\tilde s^{[r]}\colon X \to X^{[r]}$ by
induction as\footnote{For coskeletal semisimplicial complexes, one has
$\tilde s^{[r]} = \delta^{[r]} \circ s^{[0]}$, where $\delta^{[r]}
\colon X^{[0]} \to X^{[r]}$ is the diagonal embedding.}
$$
\tilde s^{[0]}(x) = s^{[0]}(x), \qquad \tilde s^{[r+1]}(x) = s^{[r+1]}
(\tilde s^{[r]}(x), \dots, \tilde s^{[r]}(x)).
$$
Since $q_i^{[2]} \circ \tilde s^{[2]} = \tilde s^{[1]}$, one has
isomorphisms $\opb{(\tilde s^{[2]})}\shl_{jk} \simeq \opb{(\tilde
s^{[1]})}\shl$.  Then, $a$ gives an isomorphism
\begin{equation}
\label{eq:LL=L}
\opb{(\tilde s^{[1]})}\shl \tens[\opb{(s^{[0]})}\ring] \opb{(\tilde
s^{[1]})}\shl \simeq \opb{(\tilde s^{[1]})}\shl.
\end{equation}
Since $\shl$ is invertible, there is an
$\ring_1\tens[\cring]\ring_0^\op$-module $\shl'$ such that $\shl
\tens[\ring_1] \shl' \simeq \ring_0$.  Applying the functor
$(\cdot) \tens[\opb{(s^{[0]})}\ring] \opb{(\tilde s^{[1]})}\shl'$ to
\eqref{eq:LL=L}, we get an isomorphism of $\opb{(s^{[0]})}\ring
\tens[\cring]\opb{(s^{[0]})}\ring^\op$-modules $\opb{(\tilde
s^{[1]})}\shl \simeq \opb{(s^{[0]})}\ring$. This proves (i).

To prove (ii), let us define the
maps $\sigma^{[r]}\colon X^{[r]} \to X^{[r+1]}$ by induction as\footnote{For coskeletal semisimplicial complexes, one has 
$\sigma^{[r]}(x) =
(x,s^{[0]}(q(x))) \in X^{[r]} \times_{X} X^{[0]} = X^{[r+1]}$, where
$q\colon X^{[r]}\to X$ is the composite of the $q^{[j]}_i$'s maps.}
$$
\sigma^{[-1]} = s^{[0]}, \qquad \sigma^{[r]}(x) = s^{[r+1]}\left(
\sigma^{[r-1]}(q^{[r]}_0(x)), \dots, \sigma^{[r-1]}(q^{[r]}_r(x)), x \right).
$$
Using the maps $\sigma^{[r]}$ one gets a functor $\opb{\sigma}\colon
\stkMod(\opb{(s^{[0]})}\ring) \to \stkMod(\twst)$, given by $\shg
\mapsto (\opb{(\sigma^{[0]})} \shl \tens[\opb q \opb{(s^{[0]})}\ring] \opb
q \shg, \opb{(\sigma^{[1]})} a)$.  This is well-defined, since
\eqref{eq:LFm} is obtained by applying $\opb{(\sigma^{[2]})}$ to
\eqref{eq:La}.  One checks that $\opb{\sigma}$ is a quasi-inverse to $\opb{\tilde s}$.
\end{proof}

\begin{proposition}
\begin{itemize}
\item[(i)]
Let $\twst$ be an $\cring$-twisting datum on $X$ satisfying \eqref{eq:locsectstrict}.  Then
$\stkMod(\twst)$ is a stack of $\cring$-twisted modules.
\item[(ii)]
Any stack of $\cring$-twisted modules on $X$ is $\cring$-equivalent to
$\stkMod(\twst)$ for some coskeletal $\cring$-twisting datum $\twst$ satisfying \eqref{eq:cosklocsectstrict}.
\end{itemize}
\end{proposition}

\begin{proof}
By definition, the maps $q^{[r]}$'s
admit local sections on $X$. Hence
part (i) follows from Lemma~\ref{le:globsect}.  As for (ii), it is enough to take $\twst$ as in \eqref{eq:twistdat0}.
\end{proof}

\begin{proposition}
\label{pr:morph}
Let $\rho\colon \twst' \to \twst$ be a refinement of coskeletal $\cring$-twisting
data on $X$ satisfying \eqref{eq:cosklocsectstrict}.  Then the functor $\opb \rho \colon \stkMod(\twst) \to
\stkMod(\twst')$ is an $\cring$-equivalence.
\end{proposition}

\begin{proof}
Let $\twst = (X^{[\bullet]} \to[q] X, \ring, \shl, a)$ and
$\twst' = (X^{\prime[\bullet]} \to[q'] X, \ring', \shl', a')$.
Proving that $\opb \rho$ is an equivalence is a local problem, and we may thus assume that $q'\colon X^{\prime[0]} \to X$ admits a global section.
Then $s = \rho^{[0]} \circ s'$ is a global section of $q\colon X^{[0]} \to X$. With the notations of Lemma~\ref{le:globsect}, one has 
$\tilde s = \rho \circ \tilde s'$.
Hence there is a diagram of
functors commuting up to an invertible transformation
$$
\xymatrix{ 
\stkMod(\twst) \ar[dr]_{\opb{\tilde s}} \ar[rr]^{\opb \rho} & &
\stkMod(\twst') \ar[dl]^{\opb{\tilde s'{}}} \\
& \stkMod(\opb{(s^{\prime [0]})}\ring'), 
}
$$
whose diagonal arrows are equivalences.
\end{proof}

\subsection{Classification of stacks of twisted 
modules}\label{sse:moritacohomology}

One may consider coskeletal $\cring$-twisting data $\twst = (X^{[\bullet]}
\to[q] X, \ring, \shl, a)$ as a kind of Cech cocycles attached to the
covering $q$, with \eqref{eq:La} playing the role of the cocycle
condition.  There is also a straightforward analogue to the notion of
coboundary, given by Morita theorem as follows.  Let $\twsu =
(X^{[\bullet]} \to[q] X, \ringi, \shm, b)$ be another coskeletal
$\cring$-twisting datum attached to the same covering $q$ as $\twst$. 
Let us say that $\twst$ and $\twsu$ differ by a coboundary if there
exist a pair $(\she,e)$ where $\she$ is an invertible
$\ring\tens[\cring]\ringi^\op$-module on $X^{[0]}$, and $e\colon
\shl\tens[\ring_1]\she_1 \isoto \she_0 \tens[\ringi_0] \shm$ is an
isomorphism of $\ring_0\tens[\cring]\ringi_1^\op$-modules on $X^{[1]}$
such that the following diagram on $X^{[2]}$ commutes
$$
\xymatrix{ \shl_{01} \tens[\ring_1] \shl_{12} \tens[\ring_2] \she_2
\ar[r]^{\id_{\shl_{01}}\tens e_{12}} \ar[d]^{a\tens\id_{\she_2}} & 
\shl_{01} \tens[\ring_1] \she_1
\tens[\ringi_1] \shm_{12} \ar[r]^{e_{01}\tens \id_{\shm_{12}}} & 
\she_0 \tens[\ringi_0]
\shm_{01} \tens[\ringi_1] \shm_{12} \ar[d]^{\id_{\she_0}\tens b} \\
\shl_{02} \tens[\ring_2] \she_2 \ar[rr]^{e_{02}} && \she_0
\tens[\ringi_0] \shm_{02} .}
$$
In this case, there is an $\cring$-equivalence
$\stkMod(\twsu)\to\stkMod(\twst)$ given by $(\shg,n) \mapsto
(\she\tens[\ringi]\shg, (\id_{\she_0} \tens n)\circ 
(e\tens\id_{\shg_1}))$.

Note that $\cring$-equivalence classes of stacks of $\cring$-twisted
modules are in one-to-one correspondence with this
``cohomology''. (The analogue correspondence appears 
in~\cite{Stevenson2000} for the case of bundle gerbes, and 
in~\cite{Breen-Messing2001} for general gerbes.) In fact, one checks 
that if $\twst =
(X^{[\bullet]} \to[q] X, \ring, \shl, a)$ and $\twsu = (Y^{[\bullet]}
\to[r] X, \ringi, \shm, b)$ are arbitrary coskeletal $\cring$-twisting data, then
$\stkMod(\twst)$ and $\stkMod(\twsu)$ are $\cring$-equivalent if and
only if $\twst$ and $\twsu$ differ by a coboundary on a
common refinement.  This means that there exist refinements of coskeletal
$\cring$-twisting data $\twst'\to\twst$ and $\twsu'\to\twsu$ such that $\twst'$
and $\twsu'$ are attached to the same covering, and differ by a
coboundary.

\subsection{Operations in terms of twisting 
data}\label{se:twistdataoperations}

Operations for stacks of twisted $\C$-modules were described in~\cite{Kashiwara-Schmid1994} using twisting data. We give here a similar description for general twisted modules.

\medskip
Let $\twst = (X^{[\bullet]} \to[q] X, \ring,\shl,a)$ be an
$\cring$-twisting datum on the topological space $X$.  Its opposite is
the $\cring$-twisting datum
\begin{equation}
\label{eq:stkinv}
\twst^\op = (X^{[\bullet]} \to[q] X, \ring^\op,\shl^{-1},a^{-1}),
\end{equation}
where $\shl^{-1} = \shHom[\ring_0](\shl,\ring_0)$, and $a^{-1}$ is the
inverse of the following chain of isomorphisms
\begin{equation*}
\begin{split}
\shl^{-1}_{02} &= \shHom[\ring_0](\shl_{02},\ring_0) \\
&\to[a] \shHom[\ring_0](\shl_{01}\tens[\ring_1]\shl_{12},\ring_0) \\
&\simeq \shHom[\ring_1](\shl_{12}, \shHom[\ring_0](\shl_{01},\ring_0))
\\
&\simeq \shHom[\ring_1](\shl_{12}, \ring_1) \tens[\ring_1]
\shHom[\ring_0](\shl_{01},\ring_0) \\
&= \shl^{-1}_{01}\tens[\ring_1^\op]\shl^{-1}_{12},
\end{split}
\end{equation*}
where in the last isomorphism holds because $\shl_{12}$ is a flat 
$\ring_1$-module locally of finite presentation.

Let $\twst' = (X^{\prime [\bullet]} \to[q'] X, \ring',\shl',a')$ be
another $\cring$-twisting datum on $X$.  Consider the semisimplicial
complex $X^{[\bullet]}\times_X X^{\prime [\bullet]} \to[p] X$, and denote
by $\pi^{[\bullet]}\colon X^{[\bullet]}\times_X X^{\prime [\bullet]} \to
X^{[\bullet]}$ and $\pi^{\prime [\bullet]}\colon X^{[\bullet]}\times_X
X^{\prime [\bullet]} \to X^{\prime [\bullet]}$ the natural maps.  If $\shf$ is
a sheaf on $X^{[r]}$ and $\shf'$ is a sheaf on $X^{\prime [r]}$, write for
short $\shf\tens[\cring] \shf' = \opb{(\pi^{[r]})}\shf \tens[\cring]
\opb{(\pi^{\prime[r]})}\shf'$ on $X^{[r]}\times_X X^{\prime [r]}$.  The
product of $\twst$ and $\twst'$ is the $\cring$-twisting datum on $X$
\begin{equation}
\label{eq:stktens}
\twst \tens[\cring] \twst' = (X^{[\bullet]}\times_X X^{\prime [\bullet]} \to[p] X,
\ring \tens[\cring] \ring', \shl \tens[\cring] \shl', a \tens[\cring]
a' ).
\end{equation}

Let $f\colon Y \to X$ be a continuous map of topological spaces. 
Consider the semisimplicial complex $Y\times_X X^{[\bullet]} \to[r]
Y$, and denote by $f^{[\bullet]} \colon Y\times_X X^{[\bullet]} \to
X^{[\bullet]}$ the natural maps.  The pull-back of $\twst$ by $f$ is
the $\opb f\cring$-twisting datum on $Y$
\begin{equation}
\label{eq:stkopb}
\opb f \twst = ( Y\times_X X^{[\bullet]} \to[r] Y,
\opb{(f^{[0]})}\ring , \opb{(f^{[1]})}\shl, \opb{(f^{[2]})}a ).
\end{equation}

One checks that, for $\twst$ and $\twst'$ satisfying \eqref{eq:cosklocsectstrict}, there are two $\cring$-equivalences and one $\opb f
\cring$-equivalence
\begin{eqnarray*}
\stkinv{\stkMod(\twst)} &\approx& \stkMod(\twst^\op), \\
\stkMod(\twst) \stktens \stkMod(\twst') &\approx&
\stkMod(\twst\tens[\cring]\twst'), \\
\stkopb f \stkMod(\twst) &\approx& \stkMod(\opb f \twst).
\end{eqnarray*}

Recall that a topological manifold $X$ is a paracompact Hausdorff 
topological space locally homeomorphic to $\R^n$.  In particular, $X$ 
is
locally compact.
In the context of twisting data, the sheaf theoretical operations of
Proposition~\ref{pr:operations} are easily described under the
assumption that $f\colon Y \to X$ is a morphism of topological
manifolds, and $X^{[\bullet]} \to[q] X$ and $X^{\prime[\bullet]}
\to[q'] X$ are semisimplicial complexes of topological manifolds with
submersive maps. (Note that this last requirement is automatically 
fulfilled for twisting data as in
\eqref{eq:twistdat0}.)

For example, let us describe the direct image functor $\oim f$.  With
the same notations as in \eqref{eq:stkopb}, consider the Cartesian
squares
$$
\xymatrix{ Y\times_X X^{[1]} \ar[r]^{r^{[1]}_i} \ar[d]^{f^{[1]}} &
Y\times_X X^{[0]} \ar[r]^r \ar[d]^{f^{[0]}} & Y \ar[d]^f \\
X^{[1]} \ar[r]^{q^{[1]}_i} & X^{[0]} \ar[r]^q & X. }
$$
If $(\shg,n)$ is an object of $\catMod(\opb f \twst)$, then $\oim f
(\shg,n) = (\oim{f^{[0]}} \shg, \oim f n)$, where $\oim f n$ is the
composite
\begin{equation*}
\begin{split}
\shl \tens[\ring_1] \opb{(q^{[1]}_1)} \oim{f^{[0]}} \shg & \simeq \shl
\tens[\ring_1] \oim{f^{[1]}} \opb{(r^{[1]}_1)} \shg \\
& \simeq \oim{f^{[1]}} (\opb{(f^{[1]})} \shl
\tens[{\opb{(f^{[1]})}\ring_1}] \opb{(r^{[1]}_1)} \shg) \\
& \isoto[n] \oim{f^{[1]}} \opb{(r^{[1]}_0)} \shg \\
& \simeq \opb{(q^{[1]}_0)} \oim{f^{[0]}} \shg.
\end{split}
\end{equation*}
Here, the first and last isomorphisms hold because the maps
$q^{[1]}_i$'s are submersive (and hence so are the
$r^{[1]}_i$'s), while the second isomorphism is due the fact that 
$\shl$ is a flat $\ring_1$-module locally of finite presentation, and 
hence locally a direct summand of a free $\ring_1$-module of finite 
rank.

\subsection{Complex powers of line bundles}\label{se:Cline}

Let us discuss the example of twisting data attached to line bundles.

\medskip
Let $X$ be a complex analytic manifold, and denote by $\OX$ its
structural sheaf.  Let $\pi\colon F \to X$ be a line bundle, let
$q\colon Y=F\setminus X \to X$ be the associated principal
$\C^{\times}$-bundle obtained by removing the zero-section, and denote
by $\shf$ the sheaf of sections of $\pi$.

As in Example~\ref{ex:twist}~(a2), consider the semisimplicial complex
where $X^{[r]}$ is the $(r+1)$-fold fibered product of $Y$.  For
$\lambda\in\C$, one has a local system on $X^{[1]}$
$$
\shl^{\lambda} = \opb p \C t^\lambda,
$$
where $p\colon X^{[1]} \to \C^{\times}$ is the map $(x,y) \mapsto
x/y$, and $\C t^\lambda\subset \O_{\C^\times}$ is the local system on
$\C^\times$ generated by $t^\lambda$.  This defines a $\CX$-twisting
datum
$$
\twst^\lambda = (X^{[\bullet]}\to[q] X, \C_Y, \shl^{\lambda}, a),
$$
where $a$ is given by $(c\, (x/y)^\lambda, d\, (y/z)^\lambda) \mapsto
cd\,(x/z)^\lambda$.

Denote by $\O_Y\twist\lambda$ the subsheaf of $\O_Y$ of
$\lambda$-homogeneous functions, i.e.\ solutions of $eu-\lambda$,
where $eu$ is the infinitesimal generator of the action of
$\C^{\times}$ on the fibers of $q$.  It is a $\opb q\OX$-module
locally constant along the fibers of $q$, and there is a natural
isomorphism $m \colon
\shl^\lambda\tens\opb{(q^{[1]}_1)}\O_Y\twist{\lambda} \to
\opb{(q^{[1]}_0)}\O_Y\twist{\lambda}$ on $X^{[1]}$ given by
$(c\,(x/y)^\lambda,\varphi(y)) \mapsto c\,\varphi(x)$.  This gives an
object
$$
\shf^{-\lambda} = (\O_Y\twist{\lambda}, m) \in \catMod(\OX;
\twst^{\lambda}).
$$
The choice of sign is due to the fact that there is an
isomorphism
$$
\shf \isoto \oim q\O_Y\twist{-1},
$$
given by $\varphi \mapsto (x \mapsto \varphi(q(x))/x )$, with inverse
$\psi \mapsto (x \mapsto \psi(x)\,x )$.

\section{Examples and applications}\label{se:local}

Giraud~\cite{Giraud1971} uses gerbes to define the second cohomology
of a sheaf of not necessarily commutative groups $\shg$,\footnote{We
will discuss in~\cite{D'Agnolo-Polesello} the linear analogue, where
$\shg$ is replaced by a not necessarily commutative $\cring$-algebra}
and if $\shg$ is abelian this provides a geometric description of the
usual cohomology group $H^2(X;\shg)$.  We consider here the case of a
sheaf of commutative local rings $\cring$, and recall how
$\cring$-equivalence classes of stacks of twisted $\cring$-modules are
in one-to-one correspondence with $H^2(X;\cring^\times)$.  We also discuss the examples of stacks of
twisted modules associated with inner forms of an $\cring$-algebra,
considering in particular the case of Azumaya algebras and TDO-rings.
As an application, we state a twisted version of an adjunction formula for sheaves and $D$-modules in the context of
Radon-type integral transforms.

\subsection{Twisted modules over commutative local rings}

Let $\cring$ be a sheaf of commutative rings.
With the terminology of Definition~\ref{de:morita}, an 
$\cring$-module is called invertible if
it is invertible as $\cring\tens[\cring]\cring^\op$-module.
Denote by
$\Pic(\cring)$ the set of isomorphism classes of invertible
$\cring$-modules, endowed
with the abelian group law given by tensor product over $\cring$. 
This is called the Picard group of $\cring$.

\begin{proposition}
Let $\stkm$ be a stack of twisted $\cring$-modules.  The group of
isomorphism classes of $\cring$-equivalences of $\stkm$ to itself is
isomorphic to $\Pic(\cring)$.
\end{proposition}

\begin{proof}
To an invertible $\cring$-module $\shl$, one associates the
$\cring$-functor $\varphi = \shl\tens[\cring](\cdot)$.  To an
$\cring$-equivalence $\varphi$ of $\stkm$ to itself, one associates
the invertible $\cring$-module $\shl = \shHom(\id_\stkm, \varphi)$.
\end{proof}

Let $\stkm$ be a stack of twisted $\cring$-modules, and denote by
$[\stkm]$ its $\cring$-equivalence class.  The multiplication $[\stkm]
[\stkm'] = [\stkm \stktens \stkm']$ is well defined, with identity
$[\stkMod(\cring)]$ and inverse $[\stkm]^{-1} = [\stkinv\stkm]$.  Let
us denote by $\TWS(\cring)$ the set of $\cring$-equivalence classes of
stacks of twisted $\cring$-modules endowed with this abelian group
structure. 

\begin{definition}
Let us say that $\cring$ is Picard good if invertible $\cring$-modules
are locally isomorphic to $\cring$ itself. 
\end{definition}

Recall that a sheaf of commutative rings $\cring$ is called local if 
for any
$U\subset X$ and any $r\in\cring(U)$ there exists an open covering
$\{V_{i}\}_{i\in I}$ of $U$ such that for any $i\in I$ either
$\cring/\cring r = 0$ or $\cring/\cring (1-r) = 0$ on $V_{i}$.
Sheaves of commutative local rings are examples of Picard good rings. 

In the rest of this section we assume that $\cring$ is Picard good.
Denote by $\cring^\times$ the multiplicative group of invertible
elements in $\cring$.

\begin{proposition}
\label{pr:H12pic}
\begin{itemize}
\item[(i)]
There is a group isomorphism
$\Pic(\cring) \simeq H^1(X;\cring^\times)$.
\item[(ii)]
There is a group isomorphism $\TWS(\cring) \simeq
H^2(X;\cring^\times)$.
\end{itemize}
\end{proposition}

Part (i) easily follows from the definition of Picard good. 
Part (ii) of the above proposition is proved as the analogue result 
for gerbes discussed e.g.\ in~\cite[\S2.7]{Breen1994}. Recall that $H^2(X;\cring^\times)$ is
calculated using hypercoverings, and coincides with Cech cohomology
if $X$ is Hausdorff paracompact.

\begin{definition}
Let $\stkm$ be a stack of twisted $\cring$-modules.  We say that
$\shf\in \stkm(X)$ is a locally free twisted $\cring$-module of
finite rank if there exists a covering $\{U_i\}_{i\in I}$ of $X$, and
$\cring|_{U_{i}}$-equivalences $\varphi_{i}\colon \stkm|_{U_{i}} \to
\stkMod(\cring|_{U_{i}})$, such that $\varphi_{i}(\shf|_{U_{i}})$ is a
locally free $\cring|_{U_{i}}$-module of finite rank.  More generally,
for an $\cring$-algebra $\ring$ we will speak of locally free
$\cring$-twisted $\ring$-modules of finite rank in 
$\catMod(\ring;\stkm)$.
\end{definition}

Note that if $\shf$ is a locally free twisted $\cring$-module of
finite rank, then for any $\cring|_U$-equivalence $\varphi\colon
\stkm|_U \to \stkMod(\cring|_U)$, $\varphi(\shf)$ is a locally free
$\cring$-module of finite rank.  Note also that the rank of $\shf$ is
a well defined locally constant function.

\begin{proposition}
\label{pr:twistedlocfree}
Let $\stkm$ be a stack of twisted $\cring$-modules.
\begin{itemize}
\item[(i)]
$\stkm$ is $\cring$-equivalent to $\stkMod(\cring)$ if and only if
$\stkm(X)$ has a locally free twisted $\cring$-module of rank $1$.
\item[(ii)]
More generally, $\stkm$ is $\cring$-equivalent to another stack of
twisted $\cring$-modules $\stkn$ if and only if
$\stkinv\stkm\stktens\stkn(X)$ has a locally free twisted
$\cring$-module of rank $1$.
\item[(iii)]
If $\stkm(X)$ has a locally free twisted $\cring$-module of rank $n$,
then $n$-fold product $\stkinv[n]\stkm = \stkm\stktens\cdots\stktens\stkm$ is $\cring$-equivalent to $\stkMod(\cring)$.
\end{itemize}
\end{proposition}

\begin{proof}
To a locally free twisted $\cring$-module $\shl$ of rank $1$ in
$\stkm(X)$ one associates the $\cring$-equivalence
$\shl\tens[\cring](\cdot)\colon \stkMod(\cring) \to \stkm$.  To an
$\cring$-equivalence $\varphi\colon \stkMod(\cring) \to \stkm$, one
associates the locally free twisted $\cring$-module of rank one
$\varphi(\cring)$.  This proves (i).  (ii) follows from (i).  As for
(iii), let $\shf\in \stkm(X)$ be a locally free twisted
$\cring$-module of rank $n$.  Then $\det \shf$ is a locally free
twisted $\cring$-module of rank $1$ in $\stkinv[n]\stkm$.
\end{proof}

\subsection{Twisting by inner forms}

Let $\cring$ be a Picard good sheaf of commutative rings, and let
$\ring$ be an $\cring$-algebra.  Denote by $\shAut[\cringalg](\ring)$
the sheaf of groups of automorphisms of $\ring$ as an
$\cring$-algebra, and by $\shInn(\ring)$ its normal subgroup of inner
automorphisms, i.e.\ the image of the adjunction morphism $\ad\colon
\ring^{\times} \to \shAut[\cringalg](\ring)$, $a\mapsto(b\mapsto
aba^{-1})$.

\begin{definition}
An $\cring$-algebra $\ringi$ is called an inner form of $\ring$ if
there exist an open covering $\{U_i\}_{i\in I}$ of $X$ and
isomorphisms $\theta_{i}\colon \ring|_{U_{i}} \isoto \ringi|_{U_{i}}$ of
$\cring$-algebras such that the automorphisms
$\theta_{j}^{-1}\circ\theta_{i}$ of $\ring|_{U_{ij}}$ are inner. 
\end{definition}

Isomorphism classes of inner forms of $\ring$ are classified by
$H^{1}(X;\shInn(\ring))$.

Assume that $\ring$ is a central $\cring$-algebra, i.e.\ that its 
center $Z(\ring)$ is equal to $\cring$.  (If $\ring$
is not central, the following discussion still holds by replacing
$\cring$ with $Z(\ring)$.)  Then the exact sequence
\begin{equation}
\label{eq:innershort}
1 \to \cring^{\times} \to \ring^{\times} \to[\ad] \shInn(\ring) \to 1
\end{equation}
induces the exact sequence of pointed sets
\begin{equation}
\label{eq:innerlong}
H^{1}(X;\ring^{\times}) \to[\gamma]
H^{1}(X;\shInn(\ring)) \to[\delta]
H^{2}(X;\cring^{\times}).
\end{equation}
If $\shl$ is a locally free $\ring$-module of rank one, then
$\gamma([\shl]) = [{\shEndo[\ring^\op](\shl^*)}]$.  If $\ringi$ is an
inner form of $\ring$, then $\delta([\ringi]) = [\stkm_\ringi]$, where
$\stkm_\ringi$ is the stack of twisted $\cring$-modules described in
the following proposition.

\begin{proposition}
\label{pr:inner}
Let $\ring$ be a central $\cring$-algebra, and $\ringi$ an inner form
of $\ring$.  Then there exists an $\cring$-equivalence
$\varphi\colon\stkMod(\ringi) \to \stkMod(\ring;\stkm_\ringi)$, where
$\stkm_\ringi$ is a stack of twisted $\cring$-modules.  Moreover,
$\varphi \simeq \shl_\ringi\tens[\ringi](\cdot)$ where $\shl_\ringi =
\varphi(\ringi)$ is a locally free $\cring$-twisted $\ring$-module of
rank one in $\catMod(\ring\tens[\cring]\ringi^\op;\stkm_\ringi)$, and
there is an isomorphism of $\cring$-algebras
$\ringi\simeq\shEndo[\ring^\op](\shl_\ringi^{*\ring})$, where
$\shl_\ringi^{*\ring} = \shHom[\ring](\shl_\ringi,\ring) \in
\catMod(\ringi\tens[\cring]\ring^\op;\stkinv{\stkm_\ringi})$.
\end{proposition}

\begin{proof}
Since $\ringi$ is an inner form of $\ring$, there exist an open
covering $\{U_{i}\}_{i\in I}$ of $X$, and isomorphisms $\theta_{i}\colon
\ring|_{U_{i}} \isoto \ringi|_{U_{i}}$ of $\cring$-algebras such that
$\theta_{j}^{-1}\circ\theta_{i}$ are inner.  Let $\varphi_i\colon
\stkMod(\ringi|_{U_{i}}) \to \stkMod(\ring|_{U_{i}})$ be the induced
$\cring|_{U_{i}}$-equivalences, denote by $\psi_i$ a quasi-inverse to $\varphi_i$, set $\varphi_{ij} =
\varphi_i|_{U_{ij}} \circ \psi_j|_{U_{ij}}$, and let
$\alpha_{ijk}\colon \varphi_{ij} \circ \varphi_{jk} \Rightarrow
\varphi_{ik}$ be the associated invertible transformations.  One
checks that $\varphi_{ij} \simeq \id_{\stkMod(\ring|_{U_{ij}})}$, so
that $\alpha_{ijk} \in \Endo(\id_{\stkMod(\ring|_{U_{ijk}})})^\times
\simeq \sect(U_{ijk}; \cring^\times)$.  By Proposition~\ref{pr:H12pic}~(ii), this is thus an
$\cring$-twisting datum defining a stack of twisted $\cring$-modules
$\stkm_\ringi$.  The equivalences $\varphi_i$ glue together, giving an
equivalence $\varphi\colon \stkMod(\ringi) \to
\stkMod(\ring;\stkm_\ringi)$.

The rest of the statement is a twisted version of Morita theorem. 
$\shl_\ringi = \varphi(\ringi)$ is a locally free $\cring$-twisted
$\ring$-module of rank one in $\catMod(\ring;\stkm_\ringi)$ which
inherits a compatible $\ringi^\op$-module structure by that of
$\ringi$ itself, and is such that
$\ringi\simeq\shEndo[\ring^\op](\shl_\ringi^{*\ring})$.
\end{proof}

\subsection{Azumaya algebras}

We consider here modules over Azumaya algebras as natural examples of
twisted $\cring$-modules.  Refer to~\cite{Giraud1971,Gabber1988} for
more details. See also~\cite{Caldararu2002,Donagi-Pantev2003}, where a twisted
version of the Fourier-Mukai transform is discussed, and~\cite{Kapustin2000}, for applications to mathematical physics.

\medskip

In this section we assume that $\cring$ is a sheaf of commutative
local rings on $X$.

\begin{definition}
An Azumaya $\cring$-algebra\footnote{The definition that we give here 
is good for the analytic topology, or for the \'etale topology. With 
this definition,
if $\ring$ is an Azumaya $\cring$-algebra, then the morphism
of $\cring$-algebras
$$
\ring\tens[\cring]\ring^\op \to \shEndo[\cring](\ring)
$$
given by $a\tens b \mapsto (c \mapsto acb)$ is an
isomorphism. 
For algebraic manifolds with the Zariski topology, it is this 
property which
is sometimes used to define Azumaya $\cring$-algebras when $\cring$ is the sheaf
of rings of regular functions. } is an $\cring$-algebra locally 
isomorphic
to the endomorphism algebra of a locally free $\cring$-module of
finite rank.  If the rank of such modules is constant and equal to
$n$, then one says that the Azumaya $\cring$-algebra has rank $n^2$.
\end{definition}

If $\shf$ is a locally free $\cring$-module of finite rank, then
$\cring$ and $\shEndo[\cring](\shf)$ are Morita equivalent. This is a 
basic example of Morita equivalence, and is proved by noticing that
$\shf$ itself is an invertible
$\cring\tens[\cring]\shEndo[\cring](\shf)^\op$-module (in fact, one 
has natural isomorphisms $\shf^*\tens[\cring]\shf\simeq
\shEndo[\cring](\shf)$ and $\shf\tens[{\shEndo[\cring](\shf)}]\shf^*
\simeq \cring$, where $\shf^*=\shHom[\cring](\shf,\cring)$).  It
follows that if $\ring$ is an Azumaya $\cring$-algebra then
$\stkMod(\ring)$ provides an example of stack of twisted
$\cring$-modules.  Moreover, the Skolem-Noether theorem (see
e.g.~\cite[Lemme~V.4.1]{Giraud1971}) asserts

\begin{proposition}
Any $\cring$-algebra automorphism of an Azumaya $\cring$-algebra is
inner.  In particular, Azumaya $\cring$-algebras of rank $n^2$ are inner forms of
the central $\cring$-algebra $M_n(\cring) =
\shEndo[\cring](\cring^n)$.
\end{proposition}

Set $GL_n(\cring) = M_n(\cring)^\times$, and $PGL_n(\cring) =
GL_n(\cring)/\cring^\times$.  Then the set of $\cring$-algebra
isomorphism classes of Azumaya $\cring$-algebras of rank $n^2$ is
isomorphic to $H^{1}(X;PGL_n(\cring))$.

\begin{proposition}
\label{pr:Azumaya}
Let $\ring$ be an Azumaya $\cring$-algebra of rank $n^{2}$.  Then
$\stkMod(\ring) \approx \stkm_\ring$ is a stack of twisted
$\cring$-modules, and there exists a locally free twisted
$\cring$-modules $\shf_\ring$ of rank $n$ in $\stkinv{\stkm_\ring}(X)
\approx \catMod(\ring^\op)$ such that $\ring \simeq
\shEndo[\cring](\shf_\ring)$ as $\cring$-algebras.
\end{proposition}

\begin{proof}
By Proposition~\ref{pr:inner} there exists a stack of twisted
$\cring$-modules $\stkm_\ring$, and an $\cring$-equivalence
$$
\varphi\colon \stkMod(\ring) \to \stkMod(M_n(\cring);\stkm_\ring).
$$
The functor $\cring^n \tens[M_n(\cring)] (\cdot)$ gives an
$\cring$-equivalence $\stkMod(M_n(\cring)) \to \stkMod(\cring)$.  By
\eqref{eq:modAstks}, this induces an $\cring$-equivalence
$$
\psi\colon \stkMod(M_n(\cring);\stkm_\ring) \to \stkm_\ring.
$$
Since $\ring$ is locally isomorphic to $M_n(\cring)$,
$\psi(\varphi(\ring))$ is locally isomorphic to $\cring^n$.  Set
$\shf_\ring = \shHom[\cring](\psi(\varphi(\ring)),\cring)$.
\end{proof}

With these notations, \eqref{eq:innershort} and \eqref{eq:innerlong}
read
$$
1 \to \cring^{\times} \to GL_{n}(\cring) \to PGL_{n}(\cring) \to 1,
$$
and
\begin{equation}
\label{eq:longAzu}
H^{1}(X;GL_{n}(\cring)) \to[\gamma_n] 
H^{1}(X;PGL_{n}(\cring)) \to[\delta_n]
H^{2}(X;\cring^{\times}),
\end{equation}
respectively.  If $\shf$ is a locally free $\cring$-module of rank
$n$, then $\gamma_n([\shf]) = [{\shEndo[\cring](\shf)}]$.  If $\ring$
is an Azumaya $\cring$-algebra of rank $n^2$, then $\delta_n([\ring])
= [\stkMod(\ring)]$.

One says that two Azumaya $\cring$-algebras $\ring$ and $\ring'$ are
equivalent if there exist two locally free (non twisted) 
$\cring$-modules of finite rank $\shf$ and $\shf'$ such that
$$
\ring\tens[\cring] \shEndo[\cring](\shf) \simeq \ring' \tens[\cring]
\shEndo[\cring](\shf').
$$

\begin{lemma}
Two Azumaya $\cring$-algebras are equivalent if and only if they are
Morita $\cring$-equivalent.
\end{lemma}

\begin{proof}
Since $\shEndo[\cring](\shf)$ and $\cring$ are Morita equivalent, so
are $\ring\tens[\cring] \shEndo[\cring](\shf)$ and $\ring$ by
\eqref{eq:modAstks}.  On the other hand, if $\ring$ and $\ring'$ are
Morita equivalent, then there is an $\cring$-equivalence
$\varphi\colon \stkMod(\cring) \to \stkMod(\ring) \stktens
\stkinv{\stkMod(\ring')} \approx
\stkMod(\ring\tens[\cring]\ring^{\prime\op})$.  Hence
$\ring\tens[\cring]\ring^{\prime\op} \simeq \shEndo[\cring](\shf)$ for
$\shf = \varphi(\cring)$.  Tensoring with
$\ring'$ we finally get an isomorphism $\ring\tens[\cring]
\shEndo[\cring](\ring') \simeq \ring' \tens[\cring]
\shEndo[\cring](\shf)$.
\end{proof}

Denote by $[\ring]$ the equivalence class of $\ring$.  The
multiplication $[\ring][\ring'] = [\ring\tens[\cring]\ring']$ is well
defined, with identity $[\cring]$ and inverse $[\ring]^{-1} =
[\ring^\op]$.  Denote by $\Br(\cring)$ the set of equivalence classes
of Azumaya $\cring$-algebras endowed with this abelian group law,
which is called Brauer group of $\cring$.  By the Skolem-Noether 
theorem one has a
group isomorphism
$$
\Br(\cring) \simeq \ilim[n]  H^{1}(X;PGL_n(\cring)).
$$
The limit of the maps $\delta_n$ in \eqref{eq:longAzu} gives a group
homomorphism
\begin{equation}
\label{eq:azu}
\delta \colon \Br(\cring) \to \TWS(\cring)
\end{equation}
which is described by $[\ring] \mapsto [\stkMod(\ring)]$.

\begin{proposition}
\label{pr:deltaBr}
The homomorphism $\delta$ is injective, and its image is contained in
the torsion part of $\TWS(\cring)$.
\end{proposition}

\begin{proof}
Let $\ring$ be an Azumaya $\cring$-algebra of rank $n^2$.  Let
$\shf_\ring$ be the locally free twisted $\cring$ module of rank $n$
in $\catMod(\ring^\op)$ of Proposition~\ref{pr:Azumaya}.  If
$\delta([\ring]) = [\stkMod(\ring)] = 0$, then $\shf_\ring$ is not
twisted, and $[\ring] = [\shEndo[\cring](\shf_\ring)] = 0$.  This
proves the injectivity.  As for the description of the image, by
Proposition~\ref{pr:twistedlocfree}~(iii) one has
$n\cdot[\stkMod(\ring)] = -[\stkinv[n]{\stkMod(\ring^\op)}] = 0$.
\end{proof}

\subsection{Twisted differential operators}

Rings of twisted differential operators (TDO-rings for short) were
introduced in a representation theoretical context
in~\cite{Beilinson-Bernstein1981,Beilinson-Bernstein1993}.  Modules
over TDO-rings provide another example of twisted modules, and we
recall here these facts following the presentation
in~\cite{Kashiwara1989} (see
also~\cite{Kashiwara-Tanisaki1996}). Since we deal with
complex analytic manifolds, as opposed to algebraic varieties, many
arguments are simpler than in loc.~cit. 

\medskip
Let $X$ be a complex analytic manifold, and denote by $\OX$ its
structural sheaf of holomorphic functions.  Recall that an $\OX$-ring
is a $\C$-algebra $\AX$ endowed with a morphism of $\C$-algebras
$\beta\colon\OX \to \AX$.  Morphisms of $\OX$-rings are morphisms of
$\C$-algebras compatible with $\beta$.

Denote by $\DX$ the sheaf of differential operators on $X$.  Recall
that $\DX$ is a simple $\OX$-ring with center $\CX$.

\begin{definition}
A TDO-ring on $X$, short for ring of twisted differential operators,
is an $\OX$-ring locally isomorphic to $\DX$ as $\OX$-ring.
\end{definition}

A TDO-ring $\AX$ has a natural increasing exhaustive filtration
defined by induction by $F_{-1}\AX = 0$, $F_{m+1}\AX = \{P\in\AX\colon
[P,a] \in F_m\AX \ \forall a\in\OX\}$, where $[P,Q] = PQ - QP$ is the
commutator.  Note that $F_{m+1}\AX = F_{1}\AX F_{m}\AX$ for $m\geq 0$,
and that the associated graded algebra $G\AX$ is naturally isomorphic
to $S_{\OX}(\Theta_X)$, the symmetric algebra of vector fields over
$\OX$.

\begin{proposition}
\label{pr:BB}
There are group isomorphisms
$$
\shAut[\OX\text{-ring}](\AX) \simeq d\OX \simeq \shInn(\AX).
$$
In particular, TDO-rings are inner forms of the central $\C$-algebra
$\DX$.
\end{proposition}

\begin{proof}
One has $\AX^\times = \OX^\times$.  Hence the short exact sequence
$$
1 \to \C_{X}^{\times} \to \OX^{\times} \to[d\log] d\OX \to 0
$$
gives a group isomorphism $\shInn(\AX) \simeq d\OX$.  This proves the
second isomorphism.  To prove the first, note that any $\OX$-ring
automorphism $\varphi$ of $\AX$ preserves the filtration.  Let $\omega
\in d\OX$, $P\in F_1\AX$, and denote by $\sigma_1(P)\in\Theta_X$ its
symbol of order one.  Then $P \mapsto P + \langle \sigma_1(P),\omega
\rangle$ extends uniquely to an $\OX$-ring automorphism of $\AX$.  On
the other hand, to an $\OX$-ring automorphism $\varphi$ of $\AX$ one
associates the closed form $\theta \mapsto \varphi(\tilde\theta) -
\tilde\theta$, where $\tilde\theta\in F_1(\AX)$ is such that
$\sigma_1(\tilde\theta) = \theta$.
\end{proof}

Let $\shf$ be a locally free $\OX$-module of rank one, and set $\shf^*
= \hom[\OX](\shf,\OX)$.  Then the basic example of TDO-ring is given
by
$$
\D_{\shf} = \shf \otens \DX \otens \shf^*,
$$
where $(s\otimes P\otimes s^*)\cdot(t\otimes Q\otimes t^*) = s \otimes
P \langle t, s^* \rangle Q \otimes t^*$.  Equivalently, $\D_{\shf}$ 
is the sheaf of differential endomorphisms of $\shf$, i.e.\ 
$\C$-endomorphisms $\varphi$ such that for any $s\in\shf$ there 
exists $P\in\DX$ with $\varphi(as) = P(a)s$ for any $a\in\OX$. More 
generally, for
$\lambda\in \C$ one has the TDO-ring
$$
\D_{\shf^{\lambda}} = \shf^{\lambda} \otens \DX \otens
\shf^{-\lambda},
$$
where $\shf^\lambda$ was described in Section~\ref{se:Cline}.  By
definition, sections of $\D_{\shf^{\lambda}}$ are locally of the form
$s^{\lambda} \tens P \tens s^{-\lambda}$, where $s$ is a nowhere
vanishing local section of $\shf$, with the gluing condition
$s^{\lambda}\tens P\tens s^{-\lambda}=t^{\lambda} \tens Q \tens
t^{-\lambda}$ if and only if $Q = (s/t)^{\lambda} P 
(s/t)^{-\lambda}$. 
This is independent from the choice of a branch for the ramified 
function $s^\lambda$.

\begin{proposition}
\label{pr:TDO}
Let $\AX$ be a TDO-ring.  Then there exists a stack of twisted
$\CX$-modules $\stkm_\AX$ such that $\stkMod(\AX)$ is $\C$-equivalent
to $\stkMod(\DX;\stkm_\AX)$.  Moreover, in
$\catMod(\OX;\stkinv{\stkm_\AX})$ there exists a locally free
$\CX$-twisted $\OX$-module of rank one $\OAX$, such that $\AX \simeq
\D_{\OAX}$ as $\OX$-rings.
\end{proposition}

\begin{proof}
By Proposition~\ref{pr:inner} there exists a stack of twisted
$\CX$-modules $\stkm_\AX$, and a $\C$-equivalence
$$
\varphi\colon \stkMod(\AX) \to \stkMod(\DX;\stkm_\AX).
$$
Since $\DX^\times = \OX^\times$, any locally free $\CX$-twisted
$\DX$-module of rank one is isomorphic to $\DX\tens[\OX]\shf$ for a
locally free $\CX$-twisted $\OX$-module of rank one $\shf$.  In
particular, $\varphi(\AX) \simeq \DX\tens[\OX]\shf_\AX$ for some
$\shf_\AX\in\catMod(\OX;\stkm_\AX)$, and we set $\OAX =
\shHom[\OX](\shf_\AX,\OX)$.
\end{proof}

\begin{conjecture}
Any stack of $\CX$-twisted $\DX$-modules is $\C$-equivalent to a stack
of the form $\stkMod(\DX;\stkm)$ for some stack of twisted 
$\CX$-modules $\stkm$. 
\end{conjecture}

Denote by $\Omega_X$ the sheaf of differential forms of top degree,
and recall that there is a natural isomorphism of $\OX$-rings $\DX^\op
\simeq \D_{\Omega_X}$.  To the TDO-rings $\AX$ and $\BX$ one
associates the TDO-rings
\begin{equation}
     \label{eq:TDOtens} \AX \tdotens \BX = \shEndo[{\AX \ctens
     \BX}](\AX \otens \BX), \qquad \tdoinv{\AX} =
     \OvX^{*}\otens\AX^{\op}\otens\OvX,
\end{equation}
where in the right-hand-side of the first equation $\AX$ and $\BX$ are
regarded as $\OX$-modules by left multiplication and $\AX \otens \BX$
is regarded as an $\AX \ctens \BX$-module by right multiplication.
Note that if $\shf$ and $\shf'$ are locally free twisted $\OX$-modules
of rank one, then
\begin{eqnarray*}
\tdoinv{\D_{\shf}} &\simeq& \D_{\shf^*}, \\
\D_{\shf} \tdotens \D_{\shf'} &\simeq& \D_{\shf\otens\shf'},
\end{eqnarray*}
where $\shf^* = \shHom[\OX](\shf,\OX)$.  

Let us denote by $[\AX]$ the $\OX$-ring isomorphism class of $\AX$.  
The
multiplication $[\AX][\BX] = [\AX \tdotens \BX]$ is well defined, with
identity $[\DX]$ and inverse $[\AX]^{-1} = [\tdoinv\AX]$.  Let us
denote by $\TDO(\OX)$ the set of $\OX$-ring isomorphism classes of
TDO-rings, endowed with this abelian group law.  As a corollary of
Proposition~\ref{pr:BB}, we get

\begin{proposition}
There is a group isomorphism $\TDO(\OX) \simeq H^{1}(X;d\OX)$.
\end{proposition}

For inner forms of $\DX$, the short exact sequence
\eqref{eq:innershort} reads
$$
1 \to \C_{X}^{\times} \to \OX^{\times} \to[d\log] d\OX \to 0.
$$
It induces the long exact cohomology sequence
$$
H^{1}(X;\C_{X}^{\times}) \to H^{1}(X;\OX^{\times}) \to H^{1}(X;d\OX)
\to H^{2}(X;\C_{X}^{\times}) \to H^{2}(X;\OX^{\times}),
$$
which may be written as the exact sequence of groups
$$
\Pic(\CX) \to \Pic(\OX) \to[\gamma] \TDO(\OX) \to[\delta] \TWS(\CX)
\to \TWS(\OX),
$$
where $\gamma([\shf]) = [\D_{\shf}]$ and $\delta([\AX]) =
[\stkm_\AX]$.  Note that \eqref{eq:TDOtens} implies the relations
$[\AX^\op] = \gamma([\Omega_X])-[\AX]$, and $[\stkm_\AX] =
-[\stkm_{\AX^\op}]$.  Note also that the complex span of the image of
$\gamma$ is described by $\lambda\cdot \gamma([\shf]) =
[\D_{\shf^\lambda}]$, for $\lambda\in\C$.

\begin{example}
Let $X = \PP$ be a complex finite dimensional projective space.  Then
the above long exact sequence reads $0\to\Z\to\C\to\C/\Z\to 0$. 
Denote by $\OP\twist{-1}$ the tautological line bundle, and for
$\lambda\in\C$ set $\OP\twist{\lambda} = (\OP\twist{-1})^{-\lambda}$. 
Then any TDO-ring on $\PP$ is of the form $\D_{\OP\twist{\lambda}}$
for some $\lambda$, and $[\stkm_{\D_{\OP\twist{\lambda}}}] =
[\stkm_{\D_{\OP\twist{\mu}}}]$ if and only if $\lambda-\mu\in\Z$.  In
this case, an equivalence $\stkMod(\D_{\OP\twist{\lambda}}) \to
\stkMod(\D_{\OP\twist{\mu}})$ is given by $\OP\twist{\mu-\lambda}
\tens[\OX] (\cdot)$.
\end{example}

\subsection{Twisted $\D$-module operations}

We recall here the twisted analogue of $D$-module operations,
following~\cite{Kashiwara1989,Kashiwara-Tanisaki1996}. (We do not recall here the classical formalism of operations for
$D$-modules, referring instead to~\cite{Kashiwara2000,Kashiwara1995}.)

\medskip
Besides the internal operations for TDO-rings recalled in
\eqref{eq:TDOtens}, there is an external operation defined as 
follows. 
Let $f\colon Y\to X$ be a morphism of complex analytic manifolds.  To
a TDO-ring $\AX$ on $X$ one associates the TDO-ring on $Y$
$$
\tdoopb f \AX = \shEndo[\opb f \AX] (\OY \tens[\opb f \OX] \opb f
\AX),
$$
where $\OY \tens[\opb f \OX] \opb f \AX$ is regarded as a right $\opb
f \AX$-module.  One has
\begin{eqnarray*}
\tdoopb f (\tdoinv{\AX}) &\simeq& \tdoinv{(\tdoopb f \AX)}, \\
\tdoopb f(\AX\tdotens\BX) &\simeq& \tdoopb f \AX \tdotens \tdoopb f
\BX.
\end{eqnarray*}
Moreover, if $\shf$ is a locally free twisted $\OX$-module
of rank one, then
$$
\tdoopb f \D_{\shf} \simeq \D_{\oopb f \shf},
$$
where $\oopb f \shf=\OY\tens_{f^{-1}\OX}f^{-1}\shf$.

Let $f\colon Y\to X$ be a morphism of complex analytic manifolds, and 
$\AX$ a TDO-ring on $X$.
Consider the transfer modules
\begin{align*}
     \AX_{Y\rightarrow X} & = \OY \tens[\opb f \OX] \opb f \AX,
     &&\text{ an $\tdoopb f \AX\tens[\CY](\opb f \AX)^\op$-module,}\\
\AX_{X\leftarrow Y} & = \opb f \AX \tens[\opb f \OX] \Ovf, &&\text{ an
$\opb f \AX \tens[\CY] (\tdoopb f \AX)^\op$-module,}
\end{align*}
where $\Ovf = \Omega_{Y} \otens \oopb f \Omega_{X}^{*}$, 
$\Omega_{X}^{*}$ denoting the dual $\hom[\OX](\Omega_{X},\OX)$ of 
$\Omega_{X}$.  Note that if
$\shf$ is a locally free twisted $\OX$-module of rank one, then
\begin{eqnarray*}
(\D_{\shf})_{Y\rightarrow X} &\simeq& \oopb f \shf \tens[\OY]
\D_{Y\rightarrow X} \tens[\opb f \OX] \opb f \shf^*, \\
(\D_{\shf})_{X\leftarrow Y} &\simeq& \opb f \shf \tens[\opb f \OX]
\D_{X\leftarrow Y} \tens[\OY] \oopb f \shf^*,
\end{eqnarray*}
where $\D_{Y\rightarrow X}$ and $\D_{X\leftarrow Y}$ are the classical
transfer bimodules.

Let $\stkm$ and $\stkm'$ be two stacks of twisted $\CX$-modules, 
$\AX$ and $\BX$ two TDO-rings on $X$. Denote by $\BDC(\AX;\stkm)$ the 
bounded derived category of $\catMod(\AX;\stkm)= 
\stkMod(\AX;\stkm)(X)$. The
usual operations for $D$-modules extend to the twisted case, yielding
the functors
\begin{align*}
\dtens & \colon \BDC(\AX;\stkm) \times \BDC(\BX;\stkm') \to \BDC(\AX
\tdotens \BX; \stkm \stktens[\C] \stkm'), \\
\dopb f & \colon \BDC(\AX; \stkm) \to \BDC(\tdoopb f \AX ; \stkopb f
\stkm), \\
\doim f & \colon \BDC(\tdoopb f \AX ; \stkopb f \stkm) \to \BDC(\AX;
\stkm),
\end{align*}
defined by $\shm \dtens \shm' = \shm \lotens \shm'$, $\dopb f \shm =
\AX_{Y\rightarrow X} \ltens[\opb f \AX] \opb f \shm$, and $\doim f
\shn = \roim f(\AX_{X\leftarrow Y} \ltens[\tdoopb f \AX] \shn)$.  The
usual formulas, like base-change or projection formula, hold. 
Moreover, all local notions like those of coherent module, of
characteristic variety, or of regular holonomic module, still make
sense.

We will also consider the functor
$$
\rhom[\AX] \colon \BDC(\AX;\stkinv\stkm)^\op \times \BDC(\AX;\stkm') \to \BDC(\stkm \stktens[\C] \stkm').
$$

\subsection{Twisted adjunction formula}

An adjunction formula for sheaves and $D$-modules in the context of
Radon-type integral transforms was established
in~\cite{D'Agnolo-Schapira1996}.  We briefly explain here how such 
formula generalizes to
the twisted case.  Note that a twisted adjunction formula for
Poisson-type integral transforms was established
in~\cite{Kashiwara-Schmid1994}, where the group action and the
topology of functional spaces are also taken into account.

\medskip

Let $X$ and $Y$ be complex analytic manifolds, $\stkm$ a stack of
twisted $\CX$-modules, $\stkn$ a stack of twisted $\CY$-modules, $\AX$
a TDO-ring on $X$, and $\AY$ a TDO-ring on $Y$.  We will use the
notations $\stkm_\AX$ and $\OAX$ from Proposition~\ref{pr:TDO}. 
Consider the natural projections
$$
X \from[\pi_1] X\times Y \to[\pi_2] Y,
$$
and set
\begin{eqnarray*}
\stkm \stketens \stkn &=& \stkopb{\pi_1}\stkm \stktens[\C]
\stkopb{\pi_2}\stkn, \\
\AX \tdoetens \AY &=& \tdoopb{\pi_1}\AX \tdotens \tdoopb{\pi_2}\AY.
\end{eqnarray*}
To $K \in \BDC( \stkm\stketens\stkinv\stkn )$ one associates the
functor
$$
K\circ(\cdot) \colon \BDC(\stkn) \to \BDC(\stkm), \quad G \mapsto
\reim{\pi_1}(K \tens \opb{\pi_2}G).
$$
To $\shk\in \BDC(\tdoinv\AX\tdoetens\AY; \stkinv\stkm\stketens\stkn)$
one associates the functor
$$
(\cdot)\dcirc \shk \colon \BDC(\AX;\stkm) \to \BDC(\AY;\stkn), \quad
\shm \mapsto \doim{\pi_2}(\dopb{\pi_1}\shm \dtens \shk).
$$
To $F\in \BDC(\stkm_\AX\stktens[\C]\stkm)$ one associates the objects of
$\BDC(\AX;\stkm)$ defined by
$$
\shc^{\natural}(F) =
\begin{cases}
    F\tens\OAX & \text{for } \natural=\omega, \\
    F\wtens\OAX & \text{for } \natural=\infty, \\
    \thom(F',\OAX) & \text{for } \natural=-\infty, \\
    \rhom(F',\OAX) & \text{for } \natural=-\omega,
\end{cases}
$$
where $F' = \rhom(F,\CX)$, and $\thom$ and $\wtens$ are the functors
of formal and temperate cohomology
of~\cite{Kashiwara1984,Kashiwara-Schapira1996}. (One checks
that the construction in~\cite{Kashiwara-Schapira1996} of the functors
of formal and temperate cohomology, starting from exact functors 
defined on the underlying real analytic manifolds, extends to 
the twisted case.)
Hence, for $\natural = \pm\infty$ we have to assume that $F$ is
$\R$-constructible.

Let $\shm \in \BDC(\AX;\stkm)$, $\shk \in \BDC(\tdoinv\AX\tdoetens\AY;
\stkinv\stkm\stketens\stkn)$, and $G \in
\BDC(\stkm_\AY\stktens[\C]\stkn)$.  Consider the solution complex $K
=\rhom[\AXY](\shk,\OAXY)$ of $\shk$, which is an object of $\BDC(
(\stkm_\AX\stktens[\C]\stkm)\stketens\stkinv{(\stkm_\AY\stktens[\C]\stkn)} )$.

\begin{theorem}
\label{th:adj}
With the above notations, assume that $\shm$ is coherent, and $\shk$
is regular holonomic, so that $K$ is $\C$-constructible.  If $\natural
= \pm\omega$, assume that $\pi_2$ is proper on $\supp(\shk)$, and that
$\chv(\shk)\inter(T^*X\times T^*_YY)$ is contained in the zero-section
of $T^*(X\times Y)$.  If $\natural = \pm\infty$, assume instead that
$G$ is $\R$-constructible.  Then, there is an isomorphism in $\BDC(\C)$
$$
\RHom[\AX](\shm,\shc^{\pm\natural}(K\circ G))\shift{d_{X}} \simeq
\RHom[\AY](\shm\dcirc\shk,\shc^{\pm\natural}(G)),
$$
where $\shift{d_{X}}$ denotes the shift by the complex dimension of
$X$.
\end{theorem}

We do not give here the proof, which follows the same lines as the one
for the non-twisted case given in~\cite{D'Agnolo-Schapira1996,Kashiwara-Schapira1996}.

\providecommand{\bysame}{\leavevmode\hbox to3em{\hrulefill}\thinspace}

\makelastpage

\end{document}